\documentclass[11pt, sec]{article}
\usepackage{graphicx}
\usepackage{amsmath}
\usepackage{amssymb}
\usepackage{amsfonts}

\makeatletter
\@addtoreset{equation}{section}   % Makes \section reset 'equation' counter.
\makeatother

\makeatletter
\@addtoreset{figure}{section}
\def\thefigure{\thesection.\@arabic\c@figure}
\makeatother

\makeatletter
\@addtoreset{table}{section}
\def\thetable{\thesection.\@arabic\c@table}
\makeatother

\begin{document}

\title{\bf  On the Siegel-Weil Theorem for Loop Groups (I)}

\author{Howard Garland and Yongchang Zhu\thanks{The second author's research is supported by Hong Kong 
 Research Grant Council earmark grant number HKUST 604507}}
\date{}

\maketitle

 \newtheorem{theorem}{Theorem}[section]
  \newtheorem{prop}[theorem]{Proposition}
  \newtheorem{cor}[theorem]{Corollary}
  \newtheorem{lemma}[theorem]{Lemma}
  \newtheorem{defn}[theorem]{Definition}
  \newtheorem{ex}[theorem]{Example}

 \newcommand{\cx}{{\bf C}}
\newcommand{\la}{\langle}
\newcommand{\ra}{\rangle}
\newcommand{\res}{{\rm Res}}
\newcommand{\expp}{{\rm exp}}
\newcommand{\Sp}{{\rm Sp}}
\newcommand{\lgp}{\widehat{\rm G}( { F} ((t)) )}
\newcommand{\svee}{\scriptsize \vee}
\newcommand{\deff}{\stackrel{\rm def}{=}}
\maketitle

\maketitle

\section{Introduction} \label{section 1}

The notion of an snt-module and our extension (1.14) of the Siegel-Weil
Theorem (also see Theorem 8.1, \S 8) grew out of our work on the
Siegel-Weil Theorem for arithmetic quotients of loop groups, which we prove
in Part II of this paper (\cite{GZ}).
 In fact (1.14) is a vital step in our proof
of the Siegel-Weil theorem in the loop case, and as far as we know, it also
seems to be a new result for automorphic forms on certain,
finite-dimensional, non-reductive groups. 
At the same time, the theory of
automorphic forms on arithmetic quotients of loop groups, specifically a
loop version of Godement's criterion (Theorem 5.3, Part II \cite{GZ}) is used to
prove the convergence theorem (Theorem 3.3) we need for the Eisenstein
series ${\rm Et}$ defined in (3.10). At the moment, even though the statement of
Theorem 3.3 only involves finite-dimensional groups, the proof using loop
groups is the only one we have!  

  To state our main result of this part, we first recall the Siegel-Weil theorem proved in \cite{W2}.

Let $F$ be a number field, ${\bf A}$ be the ring of adeles of $F$, and  $M$ be a symplectic space over $F$
 with symplectic pairing $\la , \ra$.
 Let $(V, (\, , \, ))$  be a finite dimensional vector space over $F$ with the non-degenerate, symmetric, bilinear form
 $(\, , \, )$. The space $ M \otimes_F V$ has the symplectic form given by 
 \[  \la  u_1 \otimes v_1 , u_2 \otimes v_2 \ra =  \la u_1 , u_2 \ra ( v_1 , v_2 ) .\]
The group $Sp_{M} $ and the orthogonal group $G$ of $V$ form a dual pair  
 in  the symplectic group  $Sp_{2N} $ of $ M \otimes_F V$  (where $2 N = {\rm dim}\, M  \, {\rm dim } \, V$).
Let 
  \begin{equation}\label{decom1s} 
 M = M_-   \oplus M_+  \end{equation}
 be a direct sum of Lagrangian subspaces of $M$, then 
 \[ M \otimes V =   M_- \otimes V \oplus M_+ \otimes V \]
 is a direct sum of Lagrangian subspaces of $ M \otimes V$. The Hilbert space $L^2 ( (M_- \otimes   V)_{\bf A} )$ is an irreducible 
  unitary representation of the
 metaplectic group $\widehat{Sp}_{2N} ({\bf A} )$, which is called the Weil representation.
 The  dense subspace ${\cal S} ( ( M_-\otimes  V )_{\bf A} ) )$, formed by the 
  Schwartz-Bruhat functions on $(M_- \otimes   V)_{\bf A} $, is invariant under the action of $\widehat{Sp}_{2N} ({\bf A} )$.  
 The theta functional 
 \[ \theta : {\cal S} ( ( M_- \otimes  V)_{\bf A} ) \to  {\Bbb C} \]
 given by 
   \[ \theta ( \phi ) = \sum_{r\in  M_- \otimes V } \phi ( r ) \]
 is $Sp_{2N} ( F)$-invariant.  
 
   The group $\widehat{Sp}_{2N} ({\bf A} )$ contains the commuting pair of groups $ \widehat{Sp}_{M} ( {\bf A})$ and
  $ G ( {\bf A})$, where  
   $\widehat {Sp}_{M} ( {\bf A} )$ is the preimage of  $Sp_{M}({\bf A})$ in $ \widehat{Sp}_{2N} ({\bf A} )$.
  
 For a given $\phi \in   {\cal S} ( (M_- \otimes  V)_{\bf A} )$,
  there are  two simple ways to construct a function on $  Sp_{M} ( F ) \backslash \widehat{Sp}_{M } ( {\bf A})$: one is   
    \begin{equation}\label{I}
  {\rm I} ( \phi , g  ) \deff    \int_{  G( F )  \backslash G({\bf A}) } \theta ( \pi ( g , h ) \phi ) d h  ,  \, \, \, \, \, \, 
 g \in  \widehat{Sp}_{M } ( {\bf A}), 
   \end{equation}
 where $\pi$ is the Weil representation
 and $ d h$ is the Haar measure on $G( {\bf A} )$ such that the volume of $ G( F )  \backslash G({\bf A})$ is $1$. 
   For the second way, we first  consider the function $  g \mapsto ( \pi ( g ) \phi) ( 0)$,  which is $P({\bf A})$-invariant
  (where $P$ is the Siegel parabolic subgroup of $Sp_{M}$ that fixes $M_+$).  We then form the Eisenstein series
  
 \begin{equation}\label{eisen1}
 {\rm  E}  ( \phi ,  g ) \deff  \sum_{ r \in P(F) \backslash Sp_{M} ( F ) }   ( r g \cdot \phi ) ( 0 ) . \end{equation}
  This gives a function on $Sp_{M} ( F ) \backslash \widehat{Sp}_{M } ( {\bf A})$. 
  When ${\rm dim} \, V  > {\rm dim} \, M + 2 $, both (\ref{I}) and (\ref{eisen1}) converge, and the  
  the Siegel-Weil formula asserts that the above
 two constructions are equal, i.e. 
 \begin{equation}\label{sw1}  {\rm I} ( \phi  ,  g ) =   {\rm E} ( {\phi } ,  g )  . \end{equation}
 Weil \cite{W2} proved such an identity in more generality for dual pairs constructed 
 from semisimple algebras with involutions.  Generalizations of this formula to non-convergent cases can be found in \cite{KR}.

 Put $ {\rm I} ( {\phi} )  \deff {\rm I} ( {\phi } ,  e ) $ and
 ${\rm E} ( {\phi } )  \deff {\rm E} ( {\phi } ,  e )$, we have  
  \begin{equation}\label{sw2} 
   {\rm I} ( {\phi } ) = {\rm E} ( \phi  ) 
  \end{equation}
Since $ {\rm I } ( {\phi } ,  g ) = {\rm I} ( {g \cdot \phi } ) $ and 
$ {\rm E}  ( \phi ,  g ) = {\rm E} (  {g \cdot \phi } ) $, we see that 
 (\ref{sw1}) and (\ref{sw2}) are equivalent.  
  
\

 Since the cosets  $ P(F) \backslash Sp_{M} ( F )$ are in one-to-one correspondence with the set $Gr( M )$ of Lagrangian subspaces
 of $ M$,  via the map $ P g \mapsto  M_+ g $ (we assume the symplectic group acts on $M$ from the right).
 The Eisenstein series (\ref{eisen1}) can  also be written as a summation over 
$Gr(M )$ as follows.
 Let $\pi_- : M \to M_-$ denote the projection with respect to the decomposition ({\ref{decom1s}).
   For $ r \in P(F)\backslash Sp_{M} ( F )$, let $ U = M_+ r$ . 
 Then the symplectic pairing  $ \pi_- (  U ) \times M_+ \to F $ factors through
  a non-degenerate pairing   $ \pi_- (  U ) \times M_+ / ( M_+ \cap U )  \to F $. For each $ v \in \pi_- ( U ) $, 
  let $ \tilde v \in U$ be a lifting of $v$, write $ \tilde{v} = \tilde v_- + \tilde v_+$ according the decomposition (\ref{decom1s}), then 
 the map 
 \begin{equation}\label{rh}
 \rho:   \pi_- ( U )  \to   M_+ / ( M_+ \cap U ) , \, \, \, \, v \mapsto   \tilde v_+ + M_+ \cap U 
 \end{equation}
 is well-defined.  One can prove that 
  
  \begin{equation}\label{eq1s}
        \pi ( r ) \phi ( 0)     =   \int_{ (\pi_- ( U )  \otimes V)_{\bf A} }   \psi ( \frac 12  \la x , \rho x \ra )     \phi ( x ) dx
  \deff E(  \phi , U ) ,
 \end{equation}
 where $\psi $ is the additive character of $ {\bf A} / F $ used in the definition of the Weil representation, and the measure
 on the right hand side is the Haar measure on  $ ( \pi_- ( U )  \otimes V)_{\bf A}$ normalized by the condition that
  the covolume of the lattice $\pi_- ( U ) \otimes V $ is $1$.
  So the Eisenstein series (\ref{eisen1}) can be written as 
  \[ {\rm  E} ( \phi ) = \sum_{U\in Gr( M)} E ( \phi , U ) .\]

    \

In our generalization of the Siegel-Weil formula (\ref{sw2}), we assume the symplectic  space $M$ has an additional 
 structure, which we call an snt-module.
 The groups involved are no longer reductive groups (examples: $ Sp_n ( F[ t] / ( t^k ) )$ and $G (  F[ t] / ( t^k ))$) .
 To state our generalization, we define

\begin{defn}  By a symplectic, nilpotent $t$-module
(=snt-module) $M$, we mean an $F[t]$-module which is finite-dimensional
when considered as a vector space over $F,$ and which is equipped  a
 symplectic form $\la \, , \, \ra$ such that the following conditions are satisfied
\newline (i) there exists a positive integer $N>0$ such that 
  \[  t^{N}\cdot \xi =0, \, \, \, {\rm  for \, all }\, \, \xi \in M.  \]
\newline \noindent (ii). The operator $t$ is self-dual, i.e. 
 \begin{equation}\label{1t}  \la t\xi , \eta \ra = \la \xi  ,t \eta \ra \end{equation} 
\end{defn}

 Since $t^N=0 $ on $M$, we may regard $M$ as an $F[[t]]$-module.
  We give a simple example of an snt-module. Consider $F[[t]]$-module

\begin{equation}\label{plane}
   H_{k}=F[[t]]/(t^{k})\oplus F[[t]]/(t^{k}),  
\end{equation}
with a symplectic form $\la ,\ra$ defined by the conditions
\newline \noindent  (i). Each of the two summands
is isotropic.
\newline \noindent (ii).

$\la (t^{i},0),(0,t^{j})\ra =\left\{ 
\begin{array}{c}
0,\;i+j\neq k-1 \\ 
1,\;i+j=k-1%
\end{array}%
\right. ,$ 
where $i,j=0,1,....,k-1.$

\noindent We shall prove  later (see Lemma \ref{lemma2.2}) that every snt-module is a direct sum of the above examples.
 For a given snt-module $M$, $g \in GL ( M)$ is called an snt-module automorphism if $g$ preserves both  the $F[[t]]$-module structure
 and  the symplectic structure.
We denote by 
 \[ Sp( M , t) \]
 the group of all snt-automorphisms on $M$.

  For a 
given snt-module $M$ and  a space $ V$ with non-degenerate, bilinear, symmetric form $ (\, , \, )$, then 
 the  space
\begin{equation}\label{eq1.1s}  M \otimes_F V =  M \otimes_{F[[ t]] } V[[t]] \end{equation}
 has a natural  snt-module structure, where $V[[ t]] = V \otimes_F F[[ t]] $. And
the symplectic form is defined using the first tensor product:
\[  \la  x_1 \otimes v_1 ,  x_2 \otimes v_2 \ra = \la x_1 , x_2 \ra ( v_1 , v_2 ) \]
 and the $F[[t]]$-module 
  structure is defined using the second tensor product.
    Let $Sp_{2N} $ ($2N = {\rm dim}\, M  {\rm dim} \, V $)
 denote the symplectic group
  of the symplectic space  $ M \otimes_F V$,   the group $ Sp( M , t)$ is a subgroup of $Sp_{2N} $.
 The  orthogonal group $ G ( F)$ acts on $ M \otimes_F V $ preserving the snt-module structure. But 
  a larger group $ G ( F [[ t]])$ acts on (\ref{eq1.1s}): for $ x \otimes v \in  M \otimes_{F[[ t]] } V[[t]]$, $g\in  G ( F [[ t]])$,
\[   (  x \otimes v ) \cdot g  = x \otimes (  v \cdot g ) .\]  
   It is easy to see that the action 
preserves both the symplectic structure and the $F[[t]]$-module structure, so  
 we also have a group morphism 
 \[     G ( F[[t]] ) \to Sp_{2N} . \]
We denote  $G^q ( F[[t]])$ the image.  The subgroups $ Sp( M , t)$ and $G^q ( F[[t]])$ in 
  $ Sp( M \otimes V , t ) $   obviously commute.   Our generalization of the Siegel-Weil formula is concerned with the commuting pair
($ Sp( M , t),  G^q ( F[[t]])$), which are not reductive groups  in general.  
 Let \[ Gr ( M , t )\]  denote the set of all $F[[t]]$-stable Lagrangian subspaces of $M$, so
 $ Gr ( M , t ) \subset Gr ( M ) $. 

We take a direct sum decomposition $ M= M_- \oplus M_+ $ such that  $M_- \in Gr ( M ) , M_+ \in Gr ( M , t ) $. 
  As before
  $L^2 ( (M_-  \otimes V )_{\bf A}) $ is a representation of 
  ${\widehat Sp}_{2N} ( {\bf A})$, with the theta functional 
\[ \theta :   {\cal S} ( (M_-  \otimes  V)_{\bf A}) \to {\Bbb C} , \, \, \, \phi \mapsto \theta ( \phi ) = \sum_{r\in M_- \otimes V} \phi( r) .\] 
 For a subspace $ U \in Gr ( M) $, let $ E ( \phi , U) $ be as in (\ref{eq1s}), and we define 
\begin{equation}\label{Et} 
 {\rm Et} ( \phi ) \deff  \sum_{ W \in Gr ( M , t ) }    E( \phi ,  W ).
\end{equation}
And we define 
 \begin{equation}\label{It}
 {\rm It} ( \phi ) \deff   \int_{  G^q ( F [[ t]] ) \backslash   G^q ( {\bf A} [[ t ]] )} \theta (   h \cdot  \phi  ) d h , \end{equation}
\noindent where $dh$ denotes the Haar measure on   $G^q ( {\bf A} [[ t ]] )$ such that the volume
 of $ G^q ( F [[ t]] ) \backslash   G^q ( {\bf A} [[ t ]] )$ is $1$.
 
 By Lemma \ref{lemma2.2},  $M$ is isomorphic to a direct sum of $n$ copies of $H_k$'s:
\begin{equation}\label{M}   M \cong  H_{k_1 } \oplus \dots \oplus H_{k_n} \end{equation}
assume ${\rm dim } \, V > 6 n + 2 $, and the quadratic form $( \, )$ on $V$ is $F$-anisotropic or ${\rm dim } V  - r > \frac 1 2 {\rm dim} \, M  + 1   $,
 where $r$ is the dimension of a maximal isotropic subspace of $V$,   then 
   we have the following generalization of the Siegel-Weil formula (Theorem \ref{main}, Theorem \ref{thm8.1}) :

\begin{equation}\label{tSW}
   {\rm Et} ( \phi ) = {\rm It } ( \phi ) .
\end{equation}
The condition ${\rm dim } \, V  >6 n + 2 $ is for the convergence of $ {\rm Et} ( \phi )$, and the 
 condition that $(V, ( , ))$ is $F$-anisotropic or  ${\rm dim }\, V  - r > \frac 1 2 {\rm dim  }\, M  + 1  $ 
 is for convergence of ${\rm It } ( \phi ) $.   
  This formula reduces to the classical formula (\ref{sw2}) when $k_1 = \dots = k_n = 1 $.

In general, the $Sp( M, t)$-action on $Gr(M , t )$ is not transitive, but there are only finitely many orbits, so the sum 
  $ {\rm Et } ( \phi ) = \sum E ( U, \phi )$ is a sum of Eisenstein series 
 induced from several "parabolic" subgroups (rather than only one), and each corresponds to 
 a $Sp( M, t)$-orbit in $Gr(M, t )$.  
  
In the case that $M $ is a direct sum of $n$-copies of the snt-module $H_k$, then 
  $ Sp ( M , t ) = Sp_{2n} ( F[ t] / ( t^k ) )$ and $ G^q ( F[[ t]] ) = G ( F[ t ] / ( t^k ))$.
 The formula (\ref{tSW}) means that the Siegel-Weil formula holds for symplectic and orthogonal groups
 over $F [ t] / ( t^k )$. It is implies that  the Siegel-Weil formula holds for symplectic and orthogonal groups
 over $F [ t] / ( p ( t ) )$ for arbitrary polynomial $p( t ) $.

  \

We now give an explicit example of our formula.   Let $F = {\Bbb Q}$,  and the snt - module $M$ be 
 \begin{equation}\label{1.10}   M = {\Bbb Q} [ t ] / ( t^2 ) \oplus  {\Bbb Q} / ( t^2 )  \end{equation}
with the snt-module structure given as in (\ref{plane}).   And we take a positive 
 definite even unimodular lattice $L$ of rank $N$ with the bilinear form given by $( , )$, 
and let $V= {\Bbb Q} L$. It is well-known that $N $ is divisible by $8$.
 Let $L_1 , \dots , L_g $ be the list of the positive definite even unimodular lattices of rank $N$ (up to isomorphism), let $( )_j$ denote the 
 the pairing of $L_j $, and $|{\rm Aut}_j |$ be the order of automorphism group of $L_j$. We denote  $ 1 \oplus 0 $ and 
  $ 0\oplus 1 $ in (\ref{1.10}) by $e_1 $ and $e_2 $ respectively, then $ e_1 , te_1 , e_2 , t e_2  $ is a 
  ${\Bbb Q} $-basis of $M$. Let $M_- = {\Bbb Q}  e_1 + {\Bbb Q} t e_1 $ and $M_+ ={\Bbb Q}  te_2 + {\Bbb Q} e_2 $. 
  It is clear that $M_-$ and $M_+$ are Lagrangian subspaces of $M$ and 
\[ M = M_- \oplus M_+ .\]  
   We have the Weil representation of $Sp_{4 N } ( {\bf  A} ) $ on 
 \[ {\cal S} ( ( M_- \otimes V )_{\bf A} ) .\] 
Take $\phi = \Pi_v \phi_v \in {\cal S} ( ( M_- \otimes V )_{\bf A} )$ as follows:
  for a finite place $p$ of ${\Bbb Q}$, $\phi_p $ is the characteristic function of 
   $     e_1 \otimes L_{{\bf Z}_p } +  t e_1 \otimes  L_{{\bf Z}_p } $ (where ${\bf Z}_p$ denotes the ring of $p$-adic integers and 
  $ L_{{\bf Z}_p } = L \otimes {\bf Z}_p $);
 for the real place $\infty$ of ${\Bbb Q}$, 
 \[ \phi_{\infty } ( e_1 \otimes v_1 + t e_1 \otimes v_2 ) = e^{\pi i \tau_1 ( v_1 , v_1 ) + \pi i \frac {-1} {\tau_2} ( v_2 , v_2 ) }, \]  
where $\tau_1 , \tau_2$ are complex numbers in the upper half plane.
  With the above choice of $\phi$, our new  Siegel-Weil formula (\ref{tSW}) becomes

\begin{eqnarray}\label{1.11}
     & &1+ \frac 12 \sum_{ a \geq 1 , b\in {\bf Z} : ( a, b) = 1} \, \sum_{ m ,n \in {\bf Z}: (m, n ) = 1 } 
    ( a m^2 \tau_1 + a n^2 \tau_2 + b ) ^{- \frac N 2 } \nonumber \\
  &=& C  \sum_{j=1}^g \frac 1 { | {\rm Aut}_j | }\sum_{ u , v \in L_j , u , v \, {\rm colinear } } e^{\pi i \tau_1 ( u , u ) + \pi i \tau_2 ( v , v ) }
 \end{eqnarray}
where $ (m, n) = 1$ ( $ ( a , b)=1 $) means $m , n$ (resp. $a, b$) are relatively prime, and $ u , v $ colinear means the ${\Bbb Q}$-span of 
 $u, v$ is at most $1$-dimensional. And the constant $C$ is given by (\ref{1.13}) below.

We compare (\ref{1.11}) with the classical Siegel-Weil formula

\begin{equation}\label{1.12}
     \frac 12  \sum_{ m , n\in {\bf Z}, ( m, n) = 1}  
    (  m \tau +  n  ) ^{- \frac N 2 } 
  = C  \sum_{j=1}^g \frac 1 { | {\rm Aut}_j | }\sum_{ u \in L_j  } e^{\pi i \tau  ( u , u )_j  }
 \end{equation}
 This formula expresses the $\frac N 2$-th Eisenstein series of $SL(2 , {\Bbb Z})$ as a sum of Theta series.
The constant term in $q$-expansion of both sides gives the density formula 

\begin{equation}\label{1.13}
  1 =   C  \sum_{j=1}^g \frac 1 { | {\rm Aut}_j | } ,\end{equation}
which determines $C$. 

 If we take more general test function $\phi $,  our   Siegel-Weil formula (\ref{tSW}) is 

  \begin{eqnarray}\label{1.14}
     & &1+ \frac 12 \sum_{ a \geq 1 , b\in {\bf Z} : ( a, b) = 1} \, \sum_{ m ,n \in {\bf Z}: (m, n ) = 1 } 
    ( a m^2 \tau_{11} + 2a mn \tau_{12} + a n^2 \tau_{22} + b ) ^{- \frac N 2 } \nonumber \\
  &=& C  \sum_{j=1}^g \frac 1 { | {\rm Aut}_j | }\sum_{ u , v \in L_j , u , v \, {\rm colinear } }
   e^{\pi i \left( \tau_{11}  ( u , u ) + 2 \tau_{12} ( u , v ) +   \tau_{22} ( v , v ) \right) },
 \end{eqnarray}
where $C$ is determined by (\ref{1.13}).

\

The paper is organized as follows: In \S 2 we study the structure of an snt-module $M$ (Lemma 2.2) and the structure of the corresponding group $%
Sp(M,t)$ of snt-automorphisms of $M$ (Cor. 2.6). In \S 3 we recall basic
facts about the Weil representation and study the Eisenstein series ${\rm Et}(\phi
)$ for the Weil representation associated to an snt-module. In
particular, we state the convergence theorem (Theorem 3.3) for such
Eisenstein series. As we already stated, the proof depends on the
convergence theorem for Eisenstein series on loop groups, and will be given
in Part II of this paper (\cite{GZ}). We also state a consequence, Theorem 3.4, of
Theorem 3.3 together with Prop. 1 and Prop. 2 of \cite{W2}.

In \S 4 we extend a result in Weil \cite{W2} and obtain Theorem 4.7, which
identifies certain abstract measures associated with the map $T_{W}$ (see (3.13) for the
definition of $T_{W})$ with certain gauge measures in the sense of \cite{W2}. In 
\S 5 we obtain Theorems 5.4 and 5.8. In particular, we identify the space of 
$G(F[[t]])$ - orbits in $M_{\_}\otimes V$ with the set of pairs $(W,i)$ with 
$W\in Gr(M_{\_},t)$ (the $t$ - Grassmannian), $i\in S_{t}^{2}(W)$ such that $%
U(i)_{F}$ is non-empty ( $U(i)_{F}$ is defined in \S 5, just before the
statement of Theorem 5.8).

In \S 6 we discuss $\theta $ - series and finally, we prove Theorems 7.3 and
8.1, our versions of the Siegel-Weil theorem for snt-modules.

\

\

%
%
%
%
%Section2
%

\section{The structure of symplectic, nilpotent $t$-modules.}

In this section, we prove some results about the structure of symplectic, nilpotent $t$-modules (snt-modules)
  defined in Section 1 (Definition 1.1).

\begin{lemma}
Let $M$ be an snt-module. For all $\xi \in M$,  $k\in \mathbb{Z}_{>0}$, we
have%
\begin{equation*}
\la \xi ,t^{k}\xi \ra=0.
\end{equation*}
\end{lemma}

\noindent {\it Proof. } We have for all $\xi \in M$,  $k\in \mathbb{Z}_{\geq 0},$ 
\begin{equation*}
\la \xi ,t^{k}\cdot \xi \ra =- \la t^{k}\cdot \xi ,\xi \ra,
\end{equation*}%
since $\la \cdot ,\cdot \ra $ is skew symmetric. On the other hand%
\begin{equation*}
\la \xi ,t^{k}\cdot \xi \ra =\la t^{k}\cdot \xi ,\xi \ra 
\end{equation*}%
by (\ref{1t}). Hence
\[ \la \xi ,t^{k}\cdot \xi \ra =0,\]

 \hfill $\Box$

\begin{lemma}\label{lemma2.2}
Every snt-module $M$ is isomorphic to a direct sum%
\begin{equation*}
M\cong H_{k_{1}}\oplus ....\oplus H_{k_{n}}
\end{equation*}%
 $k_{1}\geq k_{2}\geq ....\geq k_{n}$. Where $H_k$ is given as (\ref{plane}).
 Moreover, $n$ and the $k_{i}$ are uniquely determined by $M$.
\end{lemma}

\noindent {\it Proof.}
The uniqueness is clear from the theory of elementary divisors. For $\xi \in
M$, we define the order of $\xi $ to be the smallest positive integer $k,$
such that $t^{k}\cdot \xi =0.$ Pick $\xi \in M$ of maximal order $(N$ say). Then
 the
 vectors
\begin{equation*}
\xi ,t\cdot \xi ,....,t^{N-1}\cdot \xi
\end{equation*}%
are linearly independent over $F.$ To see this, we consider $F[[t]]$-submodule
  $ F[[ t ]] \xi $.  Since $ F[[t]]$ is a PID,
  $ F[[ t ]] \xi $ is isomorphic to $ F[[ t]] / ( t^k )$. It is clear that $k= N$.

 Since $\la \, , \, \ra $ is non-degenerate, we can find $\eta \in M$, so that
\[  \la t^{N-1}\cdot \xi ,\eta \ra =1 , \, \, \, \,  \la t^{j}\cdot \xi ,\eta \ra =0,\, \, \, j=0,1,...,N-2 \]
But 
then,
\[  \la \xi ,t^{N-1}\cdot \eta \ra = \la t^{N-1}\cdot \xi ,\eta \ra = 1  \]
Hence
\begin{equation*}
t^{N-1}\cdot \eta \neq 0,
\end{equation*}%
and so%
\begin{equation*}
\eta ,t\cdot \eta ,....,t^{N-1}\cdot \eta 
\end{equation*}%
are linearly independent (by the above argument) and 
\begin{equation*}
t^{N}\cdot \eta =0
\end{equation*}%
(since $N$ was assumed the maximal order of any element of $M$). But then if%
\begin{equation*}
C_{1}=F\text{-span of }\xi ,t\cdot \xi ,....,t^{N-1}\cdot \xi ,
\end{equation*}%
\begin{equation*}
C_{2}=F\text{-span of }\eta ,t\cdot \eta ,....,t^{N-1}\cdot \eta ,
\end{equation*}%
we see that%
\begin{equation*}
H=C_{1}\oplus C_{2}
\end{equation*}%
with symplectic structure given by the restriction of $\la , \ra$ on $M$ is
isomorphic to the snt-module $H_{N}$, with $C_{1},$ $%
C_{2}$ corresponding to the two direct summands $F[[t]]/(t^{N})$.  Since $\la , \ra$
restricted to $H$ is non-degenerate, $M$ decomposes as a direct sum of
$snt$-modules%
\begin{equation*}
M \cong H\oplus H^{\bot },
\end{equation*}%
and applying the induction hypothesis to $H^{\bot },$ we obtain the lemma.  \hfill $\Box$

\

For an snt-module $M$, so $M$ is in particular a symplectic space. As in Section 1,   we let $Gr(M)$ denote the set of all 
 Lagrangian subspaces, and  $Gr( M , t ) $ denote the set of Lagrangian subspaces which are stable under the action of 
 $F[[t]]$, so $Gr( M , t )\subset Gr(M)$.  We call elements in $Gr(M, t )$ $t$-Lagrangian subspaces.
We have:

\begin{lemma}
If $ U \subset M$ is a $F[[t]]$-stable, isotropic and it is not properly contained in any other
 $F[[t]]$-stable, isotropic subspaces, then $U $ is  $t$-Lagrangian, i.e. $U\in Gr( M, t )$.
\end{lemma}

\noindent {\it Proof.}
Assume $U$ is not maximal isotropic. Then there exists $v\in M$,  $v\notin U,
$ such that%
\begin{equation*}
 \la v, u \ra =0,  \, \, \text{ all } \, u\in U.
\end{equation*}%
But then for all $u\in U,$ $i\in \mathbb{Z}_{\geq 0}$%
\begin{equation*}
\la t^{i}\cdot v, u \ra =\la v,t^{i}\cdot u \ra =0,
\end{equation*}%
since $t^{i}\cdot u \in U.$

Hence, the space
$ Span_{F[[t]]}\{v, U \}$ is  $F[[t]]$-stable, isotropic and contains $M$ properly. This contradicts the maximality of $M$.
\hfill $\Box$

\

As in Section 1,  we let $Sp(M,t)$ be 
the subgroup of all $\sigma \in Sp(M)$ such that $\sigma $ is
an $F[[t]]$-module automorphism.
For a finite dimensional  $F[[t]]$-module $U$,  a subset of non-zero elements $e_1 , \dots , e_n $ is called 
a {\it quasi-basis} of $U$ if every element in $U$ can be written as an $F[[t]]$-linear combination of $e_i$'s and 
 $ a_1 e_1 + \dots + a_n e_n = 0 $ implies that all $a_i e_i = 0 $.  For example,
  for 
\[  U = F[[t]]/ ( t^{k_1} )  \oplus   \dots \oplus F[[ t]] / ( t^{k_n}) , \]
 the set (n elements) $ ( 1 , 0 , \dots , 0 ) , ( 0 , 1 , \dots, 0 ) , \dots ( 0 , \dots , 1 )$ 
  is a quasi-basis.  If $U$ has a quasi-basis consisting of $n$-elements, then 
 the  $F$-vector space $ \bar U  \deff U / t U$ is $n$-dimensional.
 If we have two such $F[[t]]$-modules $U_1 $ and $U_2$, and  $T : U_1 \to U_2$ is a morphism, then 
  $T$ induces a $F$-linear map $ \bar{T} : {\bar U}_1 \to {\bar U}_2 $.

An snt-module $M$ with decomposition as Lemma \ref{lemma2.2} is called 
{\it homogeneous}
if $k_{1}=k_{2}=....=k_{n} = k  $.   In this case the group $Sp(M,t )$ is determined by

\begin{lemma}\label{lemma2.4a}
  Let $M$ be a homogeneous snt-module as in (\ref{M}) with $k_1 = \dots = k_n = k $.
 Then  $ Sp(M,t)$ is isomorphic to $Sp_{2n} ( F[[ t]] / ( t^k ) ) $. 
  In particular,   its reduction  ${\rm mod} \, t$
defines a surjective group homomorphism%
\begin{equation*}
\pi _{0}:Sp(M,t)\rightarrow Sp_{2n}(F),
\end{equation*}%
 \end{lemma}

\noindent{\it Proof.}   We first
  consider the free $F[[t]]/ ( t^k ) $-module $ ( F[[ t]]/ ( t^k ) )^{2n} $.  It has 
  as standard $F[[t]]/ ( t^ k ) $-valued symplectic form $ \la  \, , \, \ra^{\wedge } $.
  We define an $F$-valued symplectic form $\la \, , \, \ra $ on  $ ( F[[ t]]/ ( t^k ) )^{2n}$ by 
 \[  \la  a , b \ra = {\rm coefficient \, \, of \, \, } \, \,  t^{k-1} \, \, {\rm in } \, \, \la a , b \ra^{\wedge } .\]  
With this $ \la \, , \, \ra $, $( F[[ t]]/ ( t^k ) )^{2n} $ is an snt-module, which is clearly isomorphic to 
  $M$ in the lemma.  It follows from the construction that $Sp_{2n} ( F[[t]] / ( t^k ) )$ preserves
  the snt-module structure, so we have 
 \[  Sp_{2r} ( F[[t]] / ( t^k ) ) \subset Sp ( M , t)  .\] 
For the converse inclusion is also clear.  \hfill $\Box $

\

For a homogeneous snt-module $M$ as in Lemma \ref{lemma2.4a} , $\bar M  = M / t M$ has a symplectic structure defined as follows,
 if $ \bar a , \bar b \in \bar M$, let $a , b \in M$ be their liftings, then 

\begin{equation}\label{3pair}
   \la \bar a , \bar b \ra \deff  \la  a ,  t^{k-1} b \ra .
\end{equation}

Now we turn to the case of general (possibly non-homogeneous) snt-modules 
$M$,  with direct sum decomposition as in Lemma \ref{lemma2.2}.
 For fixed $k$,  we let
\begin{equation*}
M(k)=\oplus _{k_{i}=k}H_{k_{i}},
\end{equation*}%
so%
\begin{equation}\label{3de}
M =M(l_{1})\oplus ...\oplus M(l_{s}), \, \, \, \, \, l_{1}>l_{2}>....>l_{s},  
\end{equation}
with each $M (l_{i})$ a homogeneous snt-submodule. The $M(l_{i})$'s are
mutually orthogonal with respect to $\la , \ra $.

Now if $\sigma \in Sp(M,t),$ then $\sigma $ (acting on the right, recall)
has a block decomposition with respect to (\ref{3de}),
\begin{equation}\label{3mde}
\sigma =\left[ 
\begin{array}{c}
\sigma _{1}^{1}....\sigma _{s}^{1} \\ 
........ \\ 
\sigma _{1}^{s}....\sigma _{s}^{s}%
\end{array}%
\right]   
\end{equation}%
where $\sigma _{j}^{i}:M(l_{i})\rightarrow M(l_{j}),$ and for $\xi =(\xi
_{1},....,\xi _{s})\in M$ $(\xi _{i}\in M(l_{i}))$%
\begin{equation*}
\xi \sigma =(\xi _{1},....,\xi _{s})   \left[ 
\begin{array}{c}
\sigma _{1}^{1}....\sigma _{s}^{1} \\ 
........ \\ 
\sigma _{1}^{s}....\sigma _{s}^{s}%
\end{array}%
\right]      .
\end{equation*}

 We have also the decomposition $ {\bar M}  =  M / t M $ induced from the decomposition (\ref{3de}),

 \begin{equation}\label{3de2}
{\bar M}={\bar M}(l_{1})\oplus ...\oplus {\bar M}(l_{s}),\, \, \, \, l_{1}>l_{2}>....>l_{s},  
\end{equation} 
   So $\bar \sigma : {\bar X} \to {\bar X}$ has a block decomposition:
 \begin{equation}\label{3mde}
\bar{\sigma} =\left[ 
\begin{array}{c}
\bar{\sigma}_{1}^{1}....\bar{\sigma}_{s}^{1} \\ 
........ \\ 
\bar{\sigma}_{1}^{s}....\bar{\sigma}_{s}^{s}%
\end{array}%
\right]   
\end{equation}

\begin{lemma}\label{lowertri} The matrix  $\bar{\sigma }$ is block-upper triangular, 
i.e. 
\begin{equation}\label{3tri}  \bar{\sigma}_{i}^j = 0  \, \, {\rm for}  \, \,  i < j \end{equation} 
 and the diagonal block  $\bar{\sigma}_i^i$ is in $Sp( \bar{X} ( l_i ) ) $ (recall $\bar{X} ( l_i )$
  has the symplectic structure defined by (\ref{3pair})). 
\end{lemma}

\noindent{\it Proof.}    For  $i < j $, $v_j \in M(l_j)$, 
  we have $ t^{l_j} v_j = 0 $ and $\sigma_{i}^j $ is $t$-linear, this implies that 
   $  t^i ( v_j \sigma_{i}^j ) = 0 $, then 
\begin{equation}\label{3.7}
v_j \sigma_i^j \in t^{l_i - l_j } M(l_i ) \end{equation}
  This implies that $ \bar{\sigma}_{i}^j = 0 $. 
   For $a , b \in M(l_i ) $,  we have 
\begin{equation}\label{3.8} \la t^{l_i-1} a , b \ra =  \la  t^{l_i-1}  a \sigma , b \sigma   \ra = 
    \sum_{j=1}^s \la t^{l_i-1} a \sigma_{j}^i , b \sigma_j^i \ra , \end{equation}  
 if $j < i $, by (\ref{3.7})  we have 
\[ t^{l_i-1} a \sigma_j^i \in  t^{l_i-1} t^{l_j - l_i } M(l_j)  = t^{l_j-1} M(l_j) \]
 and 
\[   b  \sigma_j^i \in   t^{l_j - l_i } M(l_j) \]
it implies that
 \[\la t^{l_i-1} a \sigma_{j}^i , b \sigma_j^i \ra  = 0 . \]
For $ j > i $, then $l_i - 1 \geq l_j$, we have
 \[  t^{l_i -1 }a \sigma_j^i \in   t^{l_i -1 } M(l_j )  = 0 .\]  
 So (\ref{3.8}) gives that 
 \[   \la t^{l_i-1} a , b \ra  =  \la t^{l_i-1} a \sigma_i^i, b \sigma_i^i \ra ,\] 
by the definition (\ref{3pair}) of the symplectic form on ${\bar M} (l_i)$, we
 prove that $\bar{\sigma}_i^i$ is a symplectic isomorphism of  ${\bar M} (l_i)$.
\hfill $\Box $

\begin{cor}\label{cor2.6}  Let $M$ be an snt-module as in (\ref{3de}), then 
$Sp(M,t)$ is the semi-direct product%
\begin{equation*}
Sp(M,t)=N \ltimes H,
\end{equation*}%
where $N$ is the unipotent radical of $Sp(M,t)$ and%
\begin{equation*}
H\cong \Pi_{i=1}^{s}Sp_{2r_{i}}(F),
\end{equation*}
where $r_i $ is the number $H_{l_i}$'s in the decomposition of $M(l_i )$. 
\end{cor}

\

Next we discuss the classification of $t$-Lagrangian subspaces.  First for the snt-module (\ref{plane}),
let $e_1 ,e_2$ denote $(1,0), ( 0, 1)$ respectively.  For each $ 0 \leq i \leq k-1 $, let $L_i $  
 be the $F$-subspace with basis
\[  t^ie_1 , t^{i+1}e_1 , \dots , t^{k-1}e_1 ,  t^{k -i}e_2 , t^{k-i+1}e_2 , \dots , t^{k-1} e_2 .\]
It is clear that $ L_i $ is a $t$-Lagrangian subspace. For an snt-module $M$ with decomposition as in Lemma \ref{lemma2.2},
 let $L_{i_j }\subset H_{k_j }$ be the subspace described as above, it is clear that 
 \begin{equation}\label{2s}  
L_{i_1 } \oplus \dots \oplus L_{i_n} \end{equation}
 is an $t$-Lagrangian subspace of $ M$.  We have 

\begin{prop} \label{prop3c} Let $M$ be a snt-module with decomposition as in Lemma \ref{lemma2.2}, then every 
 $t$-Lagrangian subspace can be tranformed by  some $g \in Sp( M , t ) $ to an $t$-Lagrangian subspace as in (\ref{2s}).
\end{prop}

This proposition will not be used later, we skip its proof.

\

%Section3

\section{The Weil representation and $t$-Eisenstein series.}

 In this section,  we first recall some basic facts about the Weil representation associated to a symplectic space over $F$,  then we
 study  the Eisenstein series ${\rm Et} ( \phi )$ (\ref{Et}) for the Weil representations associated to 
 an snt-module.

 We shall fix a non-trivial additive character $\psi : {\bf A}\to S^1$ that is trivial on $F$.
   Let $F^{2N}$ be the standard symplectic space over $F$, and $  C= C_- \oplus C_+$
  be a direct sum into Lagrangian subspaces; then
  a two-fold cover, denoted by $\widehat{Sp}_{2N} ( {\bf A})$,  of the adelic group  $Sp_{2N} ( {\bf A} )$ acts 
  on $L^2 ( C_{-, {\bf A}} )$.  The subspace $ {\cal S} ( C_{-, {\bf A} })$, formed by the Schwartz-Bruhat functions
  is invariant under this action.  We recall now the action formula.  For $g\in Sp_{2N} ( {\bf A})$, let 

\begin{equation}\label{3dec}
\left[ 
\begin{array}{cc}
\alpha _{g} & \beta _{g} \\ 
\gamma _{g} & \delta _{g}%
\end{array}%
\right]
\end{equation}
  be the block decomposition of $g$ 
with respect to the decomposition 
  \[ {\bf A}^{2N} = C_{-, {\bf A}} \oplus C_{+, {\bf A}} .\]
 In this paper, we always assume the action  of $Sp_{2N}$ on $F^{2N}$ (as well as other symplectic group actions on symplectic spaces)
    is from the right.  So $ \gamma_g $ in (\ref{3dec}) is a map from $C_+ \to C_-$.  Let $\tilde g\in  \widehat{Sp}_{2N} ( {\bf A})$
 be a lifting of $g$.  
 For
   $\phi \in  {\cal S} ( C_{-, {\bf A}} )$, 
   $ ({\tilde g}\cdot \phi ) (x)  $ equals to
  \begin{equation}\label{3action1}
 \lambda  \int_{ {\rm Im}\, \gamma_g }    S_{ g}( x + x^* )   \phi ( x \alpha_g + x^* \gamma_g ) d ( x^*\gamma_g  ), 
 \end{equation}
where $\lambda \in {\Bbb C}^*$ is a certain scalar depending only on $\tilde g$,  $d( x^*\gamma_g  )$ is a Haar measure on $ {\rm Im}\, \gamma_g  $ and 
\[ S_{g} ( x + x^*) = \psi \left( \frac 12 \la x\alpha_g , x \beta_g \ra +
          \frac 12 \la x^* \gamma_g , x^* \delta_g \ra + \la x^* \gamma_g , x \beta_g \ra \right) ;\]
 it is easy to see that $\la x^* \gamma_g , x^* \delta_g \ra$ depends only on $ x^*\gamma_g $ (not on the choice of $x^*$), therefore
 $S_{g} ( x + x^*) $ is a function of $x$ and $ x^*\gamma_g $.  
Since $\tilde g$ is unitary,   $\lambda $ can be determined up to a factor in $S^1 $.

  Let $ P $ be the subgroup of $Sp_{2N}$ that consists of elements that maps $C_+$ to itself.  An element $g $ is in $P({\bf A})$ iff 
  $ \gamma_g = 0 $.  Then there is a lifting
 $  P({\bf A}) \subset  \widehat{Sp}_{2N} ( {\bf A})$ so that for $g \in P({\bf A})$ and $\phi \in  {\cal S} ( C_{-, {\bf A}} )$,
\begin{equation}\label{para} 
    ( g\cdot \phi ) ( x ) =   | det ( \alpha_g )|_{\bf A}^{\frac 12 } \, \psi ( \frac 12 \la x\alpha_g , x \beta_g \ra  ) \phi ( x \alpha_g ), 
\end{equation}
where the factor $| det ( \alpha_g )|_{\bf A}^{\frac 12} $ guarantees the unitarity of the operator $g$.  
 There is also a lifting  
\[   Sp_{2N}( F )   \subset  \widehat{Sp}_{2N} ( {\bf A}) \]
 such that theta functional 
  \[ \theta :     {\cal S} ( C_{-, {\bf A}} ) \to {\Bbb C} \]
 given by 
 \[  \theta ( \phi ) = \sum_{r\in  C_- } \phi ( r ) \]
 is invariant under $Sp_{2N}(F)$. The action of $Sp_{2N}(F)$ is given by 
 (\ref{3action1}) with $\lambda =1$ and the Haar measure  is 
    given by the 
  condition that the covolume of  $ ({\rm Im}\,   \gamma_g ) ( F ) $ is $1$.

\

 For a given snt-module $M$ with   
  \begin{equation}\label{M3}  
  M = H_{k_1 }  \oplus \dots \oplus H_{k_n} 
\end{equation}
with $H_{k_i}$ is as in (\ref{plane}), 
 and  a space $ V$ with non-degenerate, bilinear, symmetric form $ (\, , \, )$.
Let $G$ denote the orthogonal group of $V$.
 Recall that  
\begin{equation}\label{3.1}   M \otimes_F V =  M \otimes_{F[[ t]] } V[[t]] \end{equation} 
 has a natural  snt-module structure (Section 1).
 Let $Sp_{2N} $ (where $2N = {\rm dim}\, M  {\rm dim}  \, V $)
 denote the symplectic group
  of the symplectic space  $ M \otimes_F V$.   The group $ Sp( M , t)$ is a subgroup of $Sp_{2N} $.
  The group $ G ( F [[ t]])$ acts on $ M \otimes_{F[[ t]] } V[[t]]$ in the second factor, so we have 
   a group morphism 
 \[     G ( F[[t]] ) \to Sp (  M \otimes  V , t )  . \]
 Suppose $l= {\rm max} ( k_1 , \dots , k_n )$; 
     then the image $G^q ( F[[t]])$ of the above homomorphism is isomorphic to $G( F[[t]] /  (t^l ) )$. 
 We have a commuting pair ($Sp( M, t ), G( F[[t]] /  (t^l ) )$) in $Sp ( M \otimes_F V , t )$.

\

Suppose we have a direct sum decomposition 
 \begin{equation}\label{eq3.1.1}
 M = M_- \oplus M_+ 
\end{equation}
 such that  $M_+ , M_- \in Gr ( M , t )$.
  We put 
\[   X \deff  M_- \otimes V . \]
 The space 
  $L^2 ( X_{\bf A}) $ is a representation of 
  metaplectic group  ${\widehat Sp}_{2N} ( {\bf A})$, and we have the theta functional 
\begin{equation}\label{th}
 \theta :   {\cal S} (  X_{\bf A}) \to {\Bbb C} , \, \, \, \phi \mapsto \theta ( \phi ) = \sum_{r\in X } \phi( r) .
 \end{equation}
  Recall  the Eisenstein series  (\ref{eisen1}),  (\ref{eq1s}) for $\phi \in {\cal S} ( X_{\bf A} )$ is given by 
\[ {\rm E } ( \phi ) = \sum_{ U \in Gr ( M )} E(\phi , U) =
  \sum_{U \in Gr ( M) }  \int_{ (\pi_- ( U )  \otimes V )_{\bf A} }   \psi (  \frac 12 \la x , \rho x \ra )     \phi ( x ) dx  .   \]

 Let $\pi_- : M \to M_- $ be the projection map with respect to  (\ref{eq3.1.1}).  It gives a map
 \begin{equation}\label{3a}      Gr ( M ) \to Gr ( M_- ): \, \, \, \, \,  U \mapsto \pi_- ( U) . \end{equation}
  We wish to describe the inverse image  of a given $ W\in Gr ( M_-)$. 
  Let 
   \[ W^{\bot } =\{ x \in M_+ \, | \, \la W , x \ra = 0 \}.\]
  If $ U \in Gr ( M ) $ satisfying $  \pi_- ( U ) = W$, then 
  it is easy to see that 
  \[   W^{\bot } =  M_+ \cap U .\]
  The symplectic pairing  $ W \times M_+ \to F $ 
  factors through a non-degenerate pairing 
   \[ \la \, \, \ra :   W \times M_+/ W^{\bot } \to F .\]
  We may identify $W^*$ with $  M_+/ W^{\bot }$ using this pairing.
  Recall we have the map 
  \[  \rho_U :   W \to M_+ / W^{\bot } \] 
 as defined in (\ref{rh}).  It is easy to prove that $\rho_U $ is self-dual.
 We have 

\begin{lemma}\label{lemma3a} For a given $W\in Gr ( M_- )$, the map 
  \[ U \mapsto \rho_U \] is a bijection from 
  the set of $U\in Gr ( M )$ such that $ \pi_- ( U ) = W$ to the set self-dual linear maps
   from $ W$ to $W^* =  M_+ / W^{\bot }$
\end{lemma}

\noindent{\it Proof.} It is clear that the map $ U \mapsto \rho_U $ is one-to-one. 
   If $\rho :  W \to W^* =  M_+ / W^{\bot }$ is self-dual,
  then 
 \begin{equation}\label{3b} 
 U \deff \{  w + \rho ( w ) + W^{\bot }  \, |  \, w \in W \}  \end{equation}
 is a Lagrangian subspace of $M$, and
    $\rho_U = \rho $. So the map in the lemma is also onto.
\hfill $\Box $

\

Recall the $t$-Eisenstein series defined in (\ref{Et}) is a sub-series of ${\rm E} ( \phi )$
  given by  
 \begin{equation}\label{3t}
{\rm  Et} ( \phi )  = \sum_{ U \in Gr ( X , t )} E(\phi , U)  .    
 \end{equation}

 \noindent Since $M_-$ and $M_+$ are $F[[t]]$-submodules of $M$,  
  the projection map 
  \[ \pi_- : M \to  M_- \] 
 is an $F[[t]]$-module homomorphism. 
 For each $ U \in Gr ( M , t ) $, $\pi_- ( U ) $ is an $F[[t]]$-submodule of $M_-$. 
 Denote $Gr ( M_- , t )$ the set of $F[[t]]$-submodules of $M_-$, so we have map 
 \begin{equation}\label{3.2}  P :  Gr( M , t ) \to Gr ( M_- , t ) :  \, \, \, U \mapsto \pi_- ( U ) .\end{equation}
 For a $W\in Gr ( M_- , t )$, we wish to describe the inverse image 
 \[ P_W \deff P^{-1} ( W). \]

 \noindent For any $ U \in P_W$, i.e. $ U \in Gr ( M , t ) $ and  $\pi_- ( U) = W$, we have
 the linear map 
  \[ \rho_U :  W \to  M_+ / W^{\bot }  ,\]  
 as defined in (\ref{rh}).  As before $ \rho_U$ is self-dual. 
  We now prove $\rho_U$ is $F[t]$-linear.   If $ w \in W$, 
   since $ \pi ( U ) = W$, there is $w' \in W_+$ such that $ w + w ' \in U $.
  By our definition of $\rho_U $, $\rho_U ( w ) = w' \, \, \, mod ( W^{\bot} ) $.
  Since $U$ is $t$-Lagrangian, $ t w + t w ' \in U$, this implies 
  $\rho_U ( t w ) = t w' \, \, \, mod ( W^{\bot} ) $.

 \begin{lemma}\label{lemma3.2} For each $ U\in P_W$, 
  \[  \rho_U :  W  \to  M_+ / W^{\bot } = M_+ / ( M_+ \cap U ) \]
 is $F[[t]]$-linear  and self-dual.  Conversely for each 
   $\rho :   W  \to  M_+ / W^{\bot}  $ that is $F[[t]]$-linear and self-dual, there is a 
  unique $ U\in P_W$ such that $ \rho_U  =  \rho $.  Therefore $ U \mapsto \rho_U $ is 
  a bijection from $P_W$ to $F_W$, the space of all $ \rho :  W \to   U_+ / W^{\bot}$ that is  self-dual  and  $F[[t]]$-linear .
 \end{lemma}

This lemma is an $t$-analog of Lemma \ref{lemma3a}. Suppose $ \rho $ is $F[[t]]$-linear and self-dual,
  then $U$ given as (\ref{3b})
is the unique $t$-Lagrangian subspace such that $ \rho_U = \rho $.

\

We set 
\[  {\rm Et}_W  ( \phi   ) = \sum_{ U \in P_W } E ( \phi , U ) ,\] 
By Lemma \ref{lemma3.2}, we have 
 \begin{equation}\label{3ei}
    {\rm Et}_W  ( \phi   ) = \sum_{ \rho \in F_W } \int_{ (W\otimes V)_{\bf A} }
         \phi ( x ) \psi ( \frac 12 \la x , \rho  ( x ) \ra ) d x ,
 \end{equation}
where $d x $ denotes the Haar measure on  $(W\otimes V)_{\bf A}$ such that 
  the covolume of $ W\otimes V$ is $1$.  We have 
\[   {\rm Et} ( \phi   ) = \sum_{W\in Gr ( M_- , t) }  {\rm Et}_W  ( \phi   ). \]

\

For an $F[[t]]$-submodule $W\subset M_{-},$ we let
\begin{equation*}
S_{t}^{2}(W)\subset W\otimes _{F[[t]]}W
\end{equation*}%
denote the $F[[t]]$-submodule of symmetric tensors. We define%
\begin{equation*}
T_{W}:W\otimes_F V = W\otimes_{F[[t]]} V [[t]] \rightarrow S_{t}^{2}(W)
\end{equation*}%
by%
\begin{equation}\label{eq3.2.1}
T_{W}:\sum_{i=1}^{s}w_{i}\otimes v_{i}\mapsto
\sum_{i,j=1}^{s}(v_{i},v_{j})w_{i}\otimes  w_{j},  
\end{equation}%
where $w_{i}\in W,$ $v_{i}\in V$, $i=1,....,s,$ and where $w_{i}\otimes
 w_{j}$ denotes the tensor product of $w_{i},$ $w_{j}$ in $S_{t}^{2}(W).$
   Of course $T_{W}$
 can be extended to an adelic map
\begin{equation*}
T_{W}:(W\otimes _{F}V)_{{\bf A}}\rightarrow S_{t}^{2}(W)_{{\bf A}},
\end{equation*}
and for $r\in S_{t}^{2}(W)_{\bf A},$ we set%
\begin{equation*}
\mathcal{U}_{r}=\{x\in (W\otimes _{F}V)_{{\bf A}} \, | \, T_{W}(x)=r\}.
\end{equation*}%
We now consider ${\rm Et}_W ( \phi  )$ as defined in (\ref{3ei}). In that expression we
consider $\rho \in F_{W}$. We have a pairing 

\begin{equation}\label{3c} 
  S_t^2 ( W ) \times F_W \to F , \, \, \, \, \, \,  ( \sum w_i \otimes u_i , \rho ) = 
       \sum \la w_i , \rho (u_ i ) \ra    , \end{equation}
which is clearly non-degenerate. And we have 
  \begin{equation*}
 \la w,\rho (w) \ra = ( T_{W}(w),\rho ),  \, \, \, \, \,  w\in (W \otimes V)_{\bf A} .
\end{equation*}
   We may then rewrite the right hand side of (\ref{3ei}) as 
\begin{equation}\label{3.19}
\sum_{\rho \in F_{W}  }\int_{ (W \otimes V)_{{\bf A}} }\phi ( x )
    \psi (\frac{1}{2} ( T_{W}(x),\rho ) )dx.  
\end{equation}

The following result is a corollary of convergence of  Eisenstein series on loop groups, it will be proved 
 in part II.  

\begin{theorem}\label{conv} 
 Let $M $ be an snt-module as in (\ref{M3}), where $H_k$ is as in (\ref{plane}).
 If ${\rm dim} V > 6 n + 2 $, then the series  (\ref{3t}) converges absolutely and the convergence is uniform 
 for $\phi $ varying over a compact subsets in ${\cal S} ( X_{\bf A})$.  It follows that 
 the series (\ref{3ei})= (\ref{3.19}) converges absolutely and the convergence is uniform for
 $\phi $  varying over a compact subset of ${\cal S} ( (W \otimes V)_{\bf A} )$.
\end{theorem} 

\

We can apply  Proposition 1 and Proposition 2 of \cite{W2}. 
 Using Weil's notation as in \cite{W2}: 
\begin{equation}\label{eq3.3.1}
X=(W \otimes V)_{\bf A},\, \, \, \, \,  G=S_{t}^{2}(W)_{\bf A}, \, \, \, \, \,  \Gamma =S_{t}^{2}(W),  
  \, \, \, \, \, f=T_{W}.
\end{equation}
 and  $G^* = (F_W)_{\bf A}$, $ \Gamma_* = F_W $, where
  we regard  $G^* = (F_W)_{\bf A}$ as Pontryagin dual of $G$ by the pairing
\[    S_{t}^{2}(W)_{\bf A} \times  (F_W)_{\bf A} \to S^1 , \, \, \, \, \, \, \{ a , \rho \} = \psi ( \frac 1 2 ( a , \rho ) ) . \]
 
\noindent We have 
\begin{equation*}
 F_{\phi }^*  ( g^* ) =
   \int_{ (W \otimes V)_{{\bf A}} }\phi ( x )
    \psi (\frac{1}{2}\{T_{W}(x), g^*  \})dx = \int_{X} \phi ( x ) \{ f ( x ) , g^* \} d x  . 
\end{equation*}

Theorem \ref{conv} implies that the condition of Proposition 2 in \cite{W2} is satisfied; that is 

\[  \sum_{ r^* \in \Gamma_* } |  F_{\Phi }^*  ( g^* + \gamma^*  ) | \]

\noindent converges and the convergence is uniform as $( \phi , g^*)$ varies over a compact  
 subset of ${\cal S} ( X) \times G^* $.  By using of Proposition 1 and Proposition 2 \cite{W2}, 
we obtain  

\begin{theorem}\label{thm3.5} Suppose $ {\rm dim} V > 6 n + 2 $. 
To every $r\in S_{t}^{2}(W)_{\bf A},$ there corresponds a unique positive 
measure $\mu _{r}$ on $(W\otimes V)_{\bf A}$ whose support is
contained in $\mathcal{U}_{r},$ so that for every function $\phi  $ on $
(W\otimes V)_{\bf A}$ which is continuous with compact support, the
function $F_{\phi }(r)=\int \phi d\mu _{r}$ is continuous and satisfies%
\begin{equation*}
\int F_{\phi }(r)dr=\int \phi (x)dx
\end{equation*}%
where $dr,$ $dx$ are fixed Haar measures on $S_{t}^{2}(W)_{\bf A},$ $
(W\otimes _{F}V)_{\bf A },$ respectively. Moreover, the $\mu _{r}$'s are
tempered measures and for $\phi \in {\cal S}((W\otimes V)_{\bf A}),$ $F_{\phi }$ is continuous, is an element of $
L^{1}(S_{t}^{2}(W)_{\bf A})$, 
and is the Fourier transform of the
function $F_{\phi }^{\ast }(\cdot )$ on $S_{t}^{2}(W^{\ast })_{\bf A} $ given by%
\begin{equation*}
F_{\phi }^{\ast }({\rho})=\int_{(W\otimes V)_{\bf A }}\phi (x)\psi (\frac{1}{2} ( T_{W}(x),{\rho} ) )dx. 
\end{equation*}%
Finally,%
\begin{equation}\label{etw}
 {\rm Et}_W ( \phi  ) =\sum_{r\in S_{t}^{2}(W)}\int_{(W\otimes V)_{\bf A}}%
  {\phi}d\mu _{r},
\end{equation}%
the series on the right being absolutely convergent. 
\end{theorem}

 Since the convergence of the right hand side of (\ref{etw}) is uniform  as $\phi $ varies on a compact subset
 of ${\cal S} (( W\otimes V )_{\bf A})$,  ${\rm Et}_W $ is a tempered measure on $ (W\otimes V)_{\bf A} $. 
The formula (\ref{etw}) can be restated as:

\begin{cor}\label{cor3}  We have the identity of the tempered distributions:

 \begin{equation}\label{posi}
   {\rm Et}_W = \sum_{ r \in S_t^2 ( W)} \mu_r   .
  \end{equation}
  
\end{cor}  
  
\

\

\section{An extension of Weil's abstract lemma}

In this section we study the measure $d\mu_r $ in Theorem \ref{thm3.5}. 
We begin with a statement of Proposition 1 in \cite{W2} (page 6)

\begin{lemma}\label{lemmaweil}
Let $X$ and $G$ be two locally compact, abelian groups with fixed Haar measures $dx $, $dg $ respectively.  Let
\begin{equation*}
f:X\rightarrow G
\end{equation*}%
be a continuous map such that

\begin{enumerate}
\item[(A)] For any $\Phi \in \mathcal{S}(X),$ the function $F_{\Phi }^{\ast
} $ on $G^{\ast }$ defined by%
\begin{equation}\label{6.1}
F_{\Phi }^{\ast }(g^{\ast })=\int_{X}\Phi (x) \{ f(x),g^{\ast }\}dx,  
\end{equation}
(where $dx$ $\{ \, , \, \}$ denotes the pairing
between $G$ and $G^{\ast })$ is integrable on $G^{\ast }$ and the integral $%
\int |F_{\Phi }^{\ast }|dg^{\ast }$ (where $(dg^{\ast }$ is  Haar measure on $G^{\ast
})$ dual to $dg $)  converges uniformly on every compact subset of $\mathcal{S}(X).$
\end{enumerate}

\noindent Then one can find a uniquely determined family of positive measures $\{\mu
_{g}\}_{g\in G},$ on $X,$ where support($\mu _{g})\subseteq f^{-1}(\{g\}),$
and so that for every continuous function with compact support $\Phi $ on $%
X, $ the function $F_{\Phi }$ on $G$ defined by%
\begin{equation}\label{6.2}
F_{\Phi }(g)=\int \Phi d\mu _{g},  
\end{equation}%
is continuous and satisfies 
\begin{equation}\label{6.3}
\int F_{\Phi }dg=\int \Phi dx.  
\end{equation}%
Moreover, the measures $\mu _{g}$ are tempered measures and for $\Phi \in 
\mathcal{S}(X),$ $F_{\Phi }$ is continuous, belongs to $L^{1}(G),$ satisfies
(\ref{6.3}), and is the Fourier transform of $F_{\Phi }^{\ast }.$
\end{lemma}

\

\noindent   We call $(X,G,f)$ as in Lemma \ref{lemmaweil}, {\it an admissible triple}.

\

Let $ M , M_- , M_+ , V$ be as in Section 3. And as Section 3, $W $ denotes a $F[[t]]$-submodule of $M_-$.
 For a fixed place $v$ of $F$, we set 
\[  X_v = ( W\otimes V)_{F_v } , \, \, \, \, \, \, \, \, G_v = S_t^2 ( W)_{F_v } , \]
and  let $  T_v : X_v \to G_v $ be the $F_v$-linear extension of $T_W $ defined in (\ref{eq3.2.1}).  
The dual group $G_v^*$ of $G_v$ is identified with $(F_W)_{F_v }$. 

\begin{lemma}\label{lemma6.1}  If ${\rm dim} \, V > 6 n  + 2 $, then
          above triple $(X_v, G_v , T_v)$ is an admissible triple, equivalently, 
\[  F_{\Phi }^{\ast }(g^{\ast })=\int_{X}\Phi (x) \{ f(x),g^{\ast } \}dx,   \]
  satisfies condition (A) in Lemma \ref{lemmaweil}.  
 \end{lemma}

This lemma is an analog of Proposition 5 in \cite{W2} (page 45).  We expect that  
  the condition ${\rm dim} \, V > 6 n  + 2 $ can be replaced by the weaker condition
   ${\rm dim} \, V > 6 n  + 1 $. For our purpose, the condition 
 in the lemma is enough.  

\

\noindent{\it Proof.}  For simplicity, we write $X, G, T$ for $X_v, G_v, T_v$.
Let $X_{\bf A}  = ( W\otimes V)_{{\bf A} }$,
  $G_{\bf A}  = S_t^2 ( W)_{{\bf A}}$,  and $ T_{\bf A} :   X_{\bf A} \to G_{\bf A} $ be 
  $ T_W\otimes {\bf A}$
then we have 
 \[      X_{\bf A}  = X \times  X^c , \, \, \, \, \, \, \,  G_{\bf A} = G \times  G^c \]  
where $  X^c$ is the restricted product of $  ( W\otimes V)_{F_w }$'s for $w\ne v $,
   $ G^c $ is the restricted product of $S_t^2 ( W)_{F_w}$'s for $ w \ne v $.
  And $ T_{\bf A} = T \times T_c $, where $ T_c :  X^c \to  G^c$ is defined similarly as $T$.
Let $C$ be a compact subset of ${\cal S} ( X ) $.
  We choose a  function $ \phi_0 \in {\cal S} ( X^c )$ such that
\[  \int_{X^c }  \phi_0 ( x_c ) d x_c \ne 0 .\]
 Then $ F^*_{\phi_0 } (  g_c^* )$ given by 
\[  F^*_{\phi_0 } (  g_c^* )  =   \int_{X^c }  \phi_0 ( x_c ) 
 \{  T_c ( x_c ) , g_c^* \} d  x_c ,\] 
satisfies that $   F^*_{\phi_0 } (  0 )\ne 0 $.
 Since $  F^*_{\phi_0 } (  g_c^* )$ is continuous, we have 
\[     \int_{ G^c } | F_{\phi_0 } (  g_c^* ) | d  g_c^*  = M  \ne 0 .\]

 For each function $\phi (x)  \in {\cal S} ( X)$ ,  
$ \phi ( x ) \phi_0 (  x_c )$ is in ${\cal S} ( X_{\bf A} ) $. 
Since $C$ is a compact subset of ${\cal S} ( X)$, $ C \phi_0 $ is  a compact subset of $ {\cal S} ( X_{\bf A} ) $. 
 By Theorem \ref{conv}, and the  Proposition 2 in \cite{W2}, we know that 
 $ X_{\bf A} , G_{\bf A} , T_{\bf A} $ is an admissible triple.  We consider
 the function 
    \[    \int_X \int_{X^c } \phi ( x ) \phi_0 ( x_c ) 
  \{   T(x ) , g ^* \} \{  T_c ( x_c ) , g_c^* \} d x d  x_c = F_{\phi } ( g^* ) F_{\phi _0 } (  g_c^*) .\]       
Since $  C \phi_0 $ is compact, 
   the integral 
 \[ \int_{ G^* \times  G_c^* } |  F_{\phi } ( g^* ) F_{\phi _0 } (  g_c^*) | d g^* d g_c^* \]
converges uniformly as $ \phi $ varies over $C$.  By the Fubuni theorem,
   we have 
\begin{eqnarray*}
  & & \int_{ G^* \times  G_c^* } |  F_{\phi } ( g^* ) F_{\phi _0 } (  g_c^*) | d g^* d g_c^*  \\
  & & =\int_{ G^* } |    F_{\phi } ( g^* )| d g^*  \int_{  G_c^* } |    F_{\phi_0 } (  g_c^* )| d \bar g_c^* \\
  & & =  M \int_{ G^* } |    F_{\phi } ( g^* )| d g^* .\end{eqnarray*} 
 This implies that 
\[   \int_{ G^* } |    F_{\phi } ( g^* )| d g^* \]
 converges uniformly as $\phi $ varies on $C$. \hfill $\Box $.

 \

 Since $W$ is a finite dimensional over $F$ and $ t^N W = 0 $ for $N$ large, $W$ is isomorphic to 
  \[   F [ t] / (t^{k_1} )  e_1 \oplus \dots \oplus F [ t] / ( t^{k_m} )  e_m  \]
 as a $F[[t]]$-module, where $e_1 , \dots , e_m$ is a quasi-basis of $W$.
Let  $\bar{W}$ denote the quotient  $ W / t W  $, we have 
\[\bar W \cong  F \bar e_1 \oplus \dots \oplus F \bar e_m ,\]
 where $\bar e_i$ is the projection of $e_i$. 
 Let $\bar G_v $ denote $S^2 ( \bar W_{F_v} )$, where $S^2 ( \bar W_{F_v} )$ is the subspace of the symmetric tensors in $\bar W_{F_v} \otimes \bar W_{F_v}$.
  For simplicity, we shall write $G, X, \bar G $ for $G_v , X_v , \bar G_v $.
We have $\bar T$ given  by 
\[  \bar T:  \bar X \deff  (\bar W \otimes V)_v  \to \bar G \deff S^2 ( W_v ), \, \, \, \, \, \, \,  \sum_i u_i \otimes v_i \mapsto  \sum_{i, j} ( v_i , v_j ) u_i \otimes u_j.\]  
  The condition $ {\rm dim } V > 6 n + 2 $ implies in particular ${\rm dim} V >  6 n  +2\geq 6 m +2 $; this implies that the condition for    
 Proposition 5 in \cite{W2} (page 45) is satisfied, so $( \bar X , \bar G , \bar T ) $ is an admissible triple.
The canonical map $ W \to \bar W$ induces  surjective linear maps
\[ \pi_X:   X \to \bar X , \, \, \, \, \, \, \,  \pi_G : G \to \bar G .\]
We have the commutative diagram 
\begin{equation}\label{comm}
\begin{array}{ccc}
X  & \overset{\pi _{X}}{\longrightarrow } & \bar X \\ 
\downarrow T  &  & \downarrow \bar T  \\ 
G & \overset{\pi _{G}}{\longrightarrow } & \bar G
\end{array}
\end{equation}
 Let $ f : X \to \bar G $ denote $\pi_G \circ T = \bar T \circ \pi_X $. 

\begin{lemma}\label{lemma6.2}
  For $x\in X$, the following conditions are equivalent 
\newline (1). $ T $ is submersive at $x$. 
\newline (2). $f$ is submersive at $x$.
\newline (3). $ \bar T $ is submersive at $\pi_X ( x ) $.
\end{lemma}

\noindent {\it Proof.}  Since $\pi_X $ is linear and surjective, it is submmersive at every point.
 It follows that (2) and (3) are equivalent.  Since $\pi_G$ is linear and surjective, it is submmersive  at every 
 point, it follows that (1) implies (3).  The fact that (3) implies (1) follows directly from  Lemma \ref{submersive} in Section 5.

\

\begin{lemma}\label{lemma6.3}
 The map  $f :  X \to \bar G$ satisfies the condition (A) in Lemma \ref{lemmaweil}, so 
 $( X , \bar G , f )$ is an admissible triple.
\end{lemma}

\noindent{\it Proof.} We use the following diagram to prove the lemma:

\[  \begin{array}{ccc}
   X &  \overset{\pi _{X}}{\longrightarrow } & \bar X  \\
      & f \searrow &\downarrow \bar T \\
     & & \bar G
\end{array}
\]
  \noindent Let $K$ denote the kernel of $\pi_X$. We have a map
 from ${\cal S} ( X ) \to {\cal S} ( \bar X ) $ given by 
\begin{equation}\label{6.5} 
   \Phi \mapsto \bar \Phi ( \bar x ) = \int_{K} \Phi ( k + \bar x ) d k ,
\end{equation}  
   where $k$ denotes the Haar measure on $K$.  It is clear that this map is continuous.
  Consider
\[ F_{\Phi }^{\ast }({\bar g}^{\ast })=\int_{X}\Phi (x) < f(x) , {\bar g}^* > d x \]
In the right hand side, we integrate over $K$ first, and notice that 
 \[  <   f(x) , {\bar g}^* >  = < \bar T (\bar x ),  {\bar g}^* >  \]
where $\bar x = \pi_X ( x ) $, we get 
\[ F_{\Phi }^{\ast }({\bar g}^{\ast })=\int_{\bar X}  \bar  \Phi ( \bar x ) <\bar T ( \bar x ) ,{\bar g}^{\ast }> d\bar x ,\]
   since $(\bar X , \bar G , \bar T)$ is an admissible triple, the right hand side is in
 $ L^1 ( \bar G )$. And if $\Phi $ runs through 
 a compact subset of $ {\cal S} ( X )$, then $\bar \Phi $ which is related to $\Phi $ by (\ref{6.5}) runs through a corresponding 
 compact subset of
 $ {\cal S} (\bar X )$, so the integral 
\[ \int   | F_{\Phi }^{\ast }({\bar g}^{\ast }) | d \bar{g}^* \]
converges uniformly. This proves the lemma. 
\hfill 
  $\Box $

\

By Lemma \ref{lemmaweil}, we have a family of measures $\{ \mu_{\bar g} \}_{\bar g \in \bar G}$ on $X$, with 
\newline  ${\rm support} (\mu_{\bar g} ) \subset f^{-1} ( \bar g )$,
 such that for every $\Phi \in C_c ( X)$, $ \int \Phi d \mu_{\bar g }$ is continuous function of $\bar g$ and 
  \begin{equation}\label{6.8}
     \int_{\bar G }  \left( \int \Phi d \mu_{\bar g } \right) d \bar g = \int_{X} \Phi d x   .
  \end{equation}
On the other hand, apply Lemma \ref{lemmaweil} to the admissible triple $(\bar G , \bar X , \bar T)$,  we have a family of
 measures $\{ \mu_{\bar g}^{\bar T} \}_{\bar g \in \bar G}$ on $\bar X$, with  ${\rm support} (\mu_{\bar g}^{\bar T} )  \subset \bar T^{-1} ( \bar g )$,
 and 
  \begin{equation}\label{6.9}
      \int_{\bar G }  \left( \int \bar \Phi d \mu_{\bar g }^{\bar T} \right) d \bar g = \int_{\bar X} \bar \Phi d \bar x   .
  \end{equation}
Suppose that $\Phi $ and $\bar \Phi $ are related by (\ref{6.5}) and  the Haar measures $dg , d\bar g , dk $ are compatible so that 
 the right hand sides of (\ref{6.8}) and (\ref{6.9}) are equal.  We then have

\begin{equation}\label{6.9.1}
  \int_{\bar G }  \left( \int \Phi d \mu_{\bar g } \right) d \bar g =   \int_{\bar G }  \left( \int \bar \Phi d \mu_{\bar g }^{\bar T} \right) d \bar g
 \end{equation}

  \noindent We claim that the truth of (\ref{6.9.1}) for all $\Phi \in  C_c  ( X )$ implies that 
\begin{equation}\label{6.9.2}
  \int \Phi d \mu_{\bar g } =  \int \bar \Phi d \mu_{\bar g }^{\bar T} .
\end{equation}
   To prove this,  take arbitrary $ h ( \bar g ) \in C_c ( \bar G ) $,  let $ f^* h =  h \circ f $, replace 
  $\Phi $ in  (\ref{6.9.1}) by $ f^* h \Phi $.  We get 
\[ \int_{\bar G }   h(\bar g)  \left( \int \Phi d \mu_{\bar g } \right) d \bar g =   \int_{\bar G } h(\bar g)  \left( \int \bar \Phi d \mu_{\bar g }^{\bar T} \right) d \bar g.\]
The above is true for all $h(\bar g ) \in C_c ( \bar G)$, and 
   $\int \Phi d \mu_{\bar g }$, $ \int \bar \Phi d \mu_{\bar g }^{\bar T}$ are continuous functions of $\bar g$, so we have (\ref{6.9.2}).
  We rewrite (\ref{6.9.2}) as
\begin{equation}\label{6.10} 
  \int \Phi d \mu_{\bar g } =  \int ( \int_{K}  \Phi ( k + \bar x ) d k )d \mu_{\bar g }^{\bar T} 
\end{equation}
Recall Lemma 17 \cite{W2} (page 52),
 the support of $\mu_{\bar g }^{\bar T}$ is on the $\bar T^{-1} ( \bar g )_{\rm re} $ (the regular points (= submerssive points) in $\bar T^{-1} ( \bar g )$).
 By (\ref{6.10}), the support of $\mu_{\bar g }$ is in $ \pi_{X}^{-1} \bar T^{-1 } (\bar g)_{\rm re}$, which is precisely  
 the set of the regular points in $f^{-1} ( \bar g )$ by Lemma \ref{lemma6.2}.
We have proved 

\begin{lemma}\label{lemma6.5}
  The measure $\mu_{\bar g}$ is supported on $ f^{-1} ( \bar g )_{\rm re} $, the subset of regular points of $f^{-1} ( \bar g ) $.
\end{lemma}

 \noindent  We consider the diagram

  \[  \begin{array}{ccc}
        X & & \\  
       \downarrow T &   \searrow  f & \\
      G & \overset{\pi _{G}}{\longrightarrow } & \bar G 
    \end{array}
            \]

 \noindent   For the admissible triple $(X, G , T)$, Lemma \ref{lemmaweil} implies that we have a family of measures $\mu_g $ ($g \in G $)
 supported on $T^{-1} ( g)$ such that
  for $\Phi \in C_c (X)$,  $F_{\Phi } ( g ) \deff \int \Phi d\mu_g \in C ( G ) $ , 
we have

\begin{equation}\label{6.10.1}
  \int_G  F_{\Phi} ( g )  dg = \int_X \Phi ( x ) dx 
\end{equation}
  We take a subspace of $G$ that maps isomorphically onto $\bar G$, we denote this space by $\bar G$, so we have
  the identification $ G = K \times \bar G$, where $K$ is the kernal of $\pi_G $. 
   Since  $\Phi $ has compact support, $F_{\Phi }$ has compact support, and it is continuous, so we have  
 \[  \bar g \to  \int_K  F_{\Phi } ( \bar g + k ) dk \]
 is in $ C ( \bar G )$.
  The left hand side of (\ref{6.10.1}) can be written as 
\[  \int_{\bar G} \int_K  F_{\Phi } ( \bar g + k ) dk d\bar g ,\] 
so we have 
\begin{equation}\label{6.11} 
\int_{\bar G} \int_K  F_{\Phi } ( \bar g + k ) dk d\bar g  =  \int_X \Phi ( x )  dx 
\end{equation}
 On the other hand, use the triple $( X ,\bar G , f )$, we have by (\ref{6.8} 
\begin{equation}\label{6.12} 
\int_{\bar G}  ( \int \Phi d \mu_{\bar g}  )d\bar g  =  \int_X \Phi ( x ) dx 
\end{equation}
 Comparing (\ref{6.11}) and (\ref{6.12}), we get 
\begin{equation}\label{6.12.1}
  \int_{\bar G} \int_K  F_{\Phi } ( \bar g + k ) dk d\bar g = \int_{\bar G}  ( \int \Phi d \mu_{\bar g}  )d\bar g.
\end{equation} 
  Take an arbitrary $h\in C_c ( \bar G)$, let $ f^* h = h \circ f \in C ( X)$, and replacing
 $\Phi $ in (\ref{6.12.1}) by $ f^*h  \Phi $, we get 
\[ \int_{\bar G}  h( \bar g )  \int_K  F_{\Phi } ( \bar g + k ) dk d\bar g = \int_{\bar G} h(\bar g )  ( \int \Phi d \mu_{\bar g}  )d\bar g. \]
This is true for arbitrary $h\in C_c ( \bar G )$, and both $ \int_K  F_{\Phi } ( \bar g + k ) dk$ and  $\int \Phi d \mu_{\bar g}  $ are
 continuous functions, so we have
 
\[ \int_K  F_{\Phi } ( \bar g + k ) dk = \int \Phi d \mu_{\bar g} .\]

\noindent But the measure $\mu_{\bar g}$ is supported on $f^{-1} ( \bar g )_{\rm re} $, so it is a gauge measure (see section 5 of \cite{W2} for
 the definition of "gauge" measure),
 
\[ f^{-1} ( \bar g ) =   T^{-1} ( \pi_G^{-1} ( \bar g ) ) 
     =  T_X^{-1} ( g + K ) \]
  Note that 
  \[  f^{-1} ( \bar g )_{\rm re} = \cup_{k\in K}  T^{-1} ( g + k )_{\rm re} \]
For each given $\bar g$,  $ T^{-1} ( \bar g + k)_{\rm re} $ is non-singular subvariety of $X$, we have a gauge form $d \delta_k $ on it, we have

\[ \int \Phi d \mu_{\bar g} = \int_{K}  \int_{T^{-1} ( \bar g + k )} \Phi  d \delta_k d k \]

\noindent we obtain 

\begin{equation}\label{6.16}  \int_K   F_{\Phi } ( \bar g  + k ) dk   =  \int_K  \int {T^{-1} (\bar  g + k )} \Phi ( x ) d \delta_k d k .
\end{equation}

The above holds for arbitrary $\Phi \in C_c ( X)$, use the same method we used to deduce (\ref{6.9.2})  from (\ref{6.9.1}), 
 we deduce from (\ref{6.16}) that

 \[  F_{\Phi } ( \bar g  + k ) = \int {T^{-1}( \bar g + k )_{\rm re}} \Phi  d \delta_k  \]
So we have proved

\begin{lemma}\label{lemma6.6}
  The measure $\mu_{g, v} $ for the triple $(X_v, G_v , T_v )$ is supported on
 $T_v^{-1} ( g)_{\rm re}$ and is the gauge measure. 
\end{lemma}

Lets recall the meaning of "gauge" measure ( \cite{W2}, section 5). In the situation as Lemma \ref{lemma6.6}. 
We first take an invariant top form $\eta $ on $G$ and an invariant top form $\omega $ on $X$. 
 Let $X'$ be the open set of $X$ that consists of all the  points where $T$ is submersive.  Near each point $x \in X'$,
 there is a  form $\theta_x $ such that $ \theta \wedge T^* \eta = \omega $. For each $y\in G$, the local forms $\theta $, restrict
 to $ T^{-1} ( y )_{\rm re} = T^{-1} ( y ) \cap X'$,  give a top form $\theta_y $ on   $T^{-1} ( y )_{\rm re}$, which defines a measure
 which is equal to $\mu $ in Lemma \ref{lemma6.6}.

\

Now we consider the global situation.  Let $ X , G , T_W $ as in (\ref{eq3.3.1}). 
For each $r\in S^2_t ( W)$, we consider the inverse image $T_W^{-1} ( r )$.  Notice that $T$ is submersive at a generic point,
the space $X'$ formed by the points at which $f$ is submersive is  $F$-open in $X$. We take a top form $\eta$ over $G$ and a top form
 $\omega $ on $X$, we assume that the Tamagawa measures on $G ({\bf A})$ ($X_{\bf A}$ resp.)  with respect to $\eta $ (resp $\omega $ )
 are the Haar measure normalized by the condition that the covolume of $G(F)$ ($X(F)$) is $1$. The space $X'$ can be covered by $F$-open subsets
 $U_{\lambda }$ such that, there is a form
 $\theta_{\lambda } $  rational over $F$ satisfying
 $ \theta \wedge f^*  = \omega $. For each $i \in G(F)$, 
 then the $\theta$'s restrict on $f^{-1} (i) \cap X'$ to get a top form  $\theta_i $ on $ f^{-1} ( i) \cap X'$.  By Lemma \ref{lemma6.6},
 $ d \mu_{i, v}  $ is given by $ |\theta_i |_v$.
 Using  a similar argument as in \cite{W2}, Section 42, we can prove    
  that $1$ is a system of convergence factor of $|\theta_i|_{\bf A}$. And we have

\begin{theorem}\label{thm4.7} For each $ r \in S_t^2 ( W)$,  
 the measure $\mu_r $ in (\ref{etw}) is supported in $T_W^{-1} ( r)_{\rm re}$ and 
  it is the same as $|\theta_r |_{\bf A}$,  the measure define by the form $\theta_r $.
\end{theorem}

\

\

\section{Classification of orbits of orthogonal groups}

As in Section 3, we denote $M$ an snt-module with decomposition  $ M = M_- \oplus M_+$ into
  $t$-Lagrangian subspaces
and $V$ a finite dimensional vector space over $F$ with a non-degenerate, bilinear symmetric form $( , )$.
The orthogonal group $G ( F [[ t]] )$ acts on $ M \otimes V$, leaving the subspace $M_-\otimes V$ invariant.
  The purpose of this section to give 
 a complete set of invariants of $G(F[[t]])$-orbits in $M_-\otimes V$ (see Theorem \ref{theorem5.1} below).

\begin{defn}\label{def5.1} Let $W$ be a finitely generated  $F[[t]]$-module. 
 A submodule $L\subset W$ is called  a  primitive submodule if
 one of the following equivalent conditions is satisfied:
\newline (1). there is a complement $F[[t]]$-submodule $L'$, i.e. $W = L\oplus L'$.
\newline (2). the natural map $ L/ t L \to  W /t W$ induced form the embedding 
  $ L \hookrightarrow M$ is injective. 
\end{defn}

 \

\noindent {\it Examples.} (1) Let  $ W = F[[ t]] / ( t^{k_1} ) \oplus \dots \oplus F[[ t]] / ( t^{k_m} ) $,
   For any $ l \leq    m $,  $L =  F[[ t]] / ( t^{k_1} ) \oplus \dots \oplus F[[ t]] / ( t^{k_l} )$
  is a primitive submodule of $M$.
\newline (2).  Let $F[[t]]^m$ be a free $F[[t]]$-module of rank $m$, then 
  $F[[t]]^l$ (consists of elements with last $(m-l)$-components $0$ is a primitive submodule.
 (3). $\{0\}$ is a primitive submodule for any $M$.

\

Let $W$ be a finitely generated $F[[t]]$-module,  $e_1 , \dots ,e_m$ is called 
 a {\it quasi-basis} of $W$ if every $e_i\ne 0 $,  every element $x\in W$ can be written as 
 a $F[[t]]$-linear combination of $e_1 , \dots , e_m $, and $a_1 e_1 + \dots + a_m e_m = 0 $
 ($a_i\in F[[t]]$) implies that all $a_ie_i = 0 $. If $W$ is a finite dimensional $F[[t]]$-module, the "quasi-basis" defined above is the 
 same notion as  defined in Section 2.
 If $W$ is a free $F[[t]]$-module,
  then a quasi-basis of $W$ is the same as a basis of $W$. In the example (1) above,
  there is a quasi-basis of $m$ elements. It is clear that $e_1 , \dots , e_m $ is a quasi-basis of $W$ if and only if 
  the images of $e_1 , \dots , e_m $ in $W / t W $ form an $F$-basis of vector space $W /t W$. Therefore any 
  two quasi-bases have the same number of elements. The cardinality of a quasi-basis is called the {\it rank} 
  of $W$.

  \

Every 
 \[ x = \sum_{i} u_i \otimes v_i \in  M_- \otimes V = M_- \otimes_{F[[t]]} V[[t]]\]
 gives rise to an $F[[t]]$-linear map
\begin{equation}\label{5.1}
   f_x :    V[[t]] \to  M_-, \, \, \, \, \, \,  f_x (  v ) = \sum_i  ( v_i , v ) u_i  .
\end{equation}
  We denote by  ${\rm Im} \, f_x  $ the image of $f_x$, which is an $F[[t]]$-submodule of $M_-$.
Let $e_1 , \dots , e_m$ be a quasi-basis of ${\rm Im} \, f_x  $, then  
 \begin{equation}\label{5.1.1}
  {\rm Im } \, f_x \cong F[[t]]/(t^{k_1}) e_1 \oplus \dots \oplus F[[ t]] / ( t^{k_m} ) e_m ,
\end{equation} 
where $k_i$ is the smallest positive integer such that $t^{k_i} e_i = 0 $.

\begin{lemma}\label{lemma5.1} Let $e_1 , \dots , e_m $ be a quasi-basis of ${\rm Im}\, f_x $, and suppose
 $v_1 , \dots , v_m \in V[[t]] $ satisfy $ f_x ( v_i) =e_i $ ($i=1, \dots , m)$, then 
 ${\rm Span}_{F[[t]]}\{ v_1 , \dots ,v_m\}$ is a primitive submodule of $V[[t]]$ and 
 $v_1 , \dots ,v_m$ is a basis of ${\rm Span}_{F[[t]]}\{ v_1 , \dots ,v_m\}$.
  \end{lemma}

\noindent{\it Proof.}  Set $ L = {\rm Span}_{F[[t]]}\{ v_1 , \dots ,v_m\}$, then
  $f_x |_L : L \to {\rm Im} \, f_x $ induces a linear isomorphism:
\[ \bar {f_x} :    L / t L \to {\rm Im } \, f_x / t  ({\rm Im } \, f_x ). \]   
This implies in particular, ${\rm dim }  \, L/tL = m$, so the map $  L/tL  \to V[[ t ]] / t V[[ t]] $ induced from 
  $ L \subset V[[t]] $  is injective, so 
 $L$ is primitive submodule of $V[[t]]$.  The other conclusions are clear. 
\hfill $\Box$

\

The bilinear form $( , )$ on $V$ can be extended to a $F[[t]]$-valued bilinear form on $V[[t]] = V\otimes_F F [[ t]]$.  
This  $F[[t]]$-valued bilinear form on $V[[t]]$ is non-degenerate. It is easy to prove
 
\begin{lemma}\label{lemma5.2} Let $e_1 , \dots , e_m $ be a quasi-basis of ${\rm Im}\, f_x $, then 
 there are elements $ w_1 , \dots , w_m \in V[[t]]$ such that
\newline (1) ${\rm Span}_{F[[t]]}\{ w_1 , \dots ,w_m\}$ is a primitive submodule of $V[[t]]$ and 
 $w_1 , \dots ,w_m$ is basis of ${\rm Span}_{F[[t]]}\{ w_1 , \dots ,w_m\}$.
\newline (2) \[ x = e_1 \otimes w_1 + \dots + e_m \otimes w_m .\]
  \end{lemma}

\noindent {\it Proof.} Choose $v_i \in V[[ t]] $ ( $i=1 , \dots , m $) such that $ f_x ( v_i ) = e_ i $.
 By Lemma \ref{lemma5.1},  $ {\rm Span}_{F[[t]]} \{ v_1 , \dots , v_m \} $ is a primitive submodule of
 $ V[[t]]$ and $v_1 , \dots , v_m $ is a basis.  It is clear that 
  \[ V[[ t ]] = {\rm Span}_{F[[t]]} \{ v_1 , \dots , v_m \} \oplus ker ( f_x )  .\] 
Let $ v_{m+1} , \dots , v_{N } $ be a basis of $ ker ( f_x ) $. Then 
  $ v_1 , \dots , v_N $ is basis of $V[[t]]$. Now let $w_1 , \dots , w_N $ be the dual basis if
  $ v_1 , \dots , v_N$, i.e., $ ( v_i , w_j ) = \delta_{ij}$. 
This is clear that $ v = \sum_{i=1}^m e_i \otimes w_i $.
 \hfill $\Box $

From Lemma \ref{lemma5.2}, we see that $x\in {\rm Im}\, f_x \otimes_{F[[t]]} V[[ t]] $. 
  We define a map
\[ T:  {\rm Im} \, f_x  \otimes_{F[[t]]} V [[ t]] \to  S_t^2 ( {\rm Im } \, f_x ) , \, \, \, \, \, \,  \sum_i^N  u_i \otimes v_i \mapsto 
\sum_{i, j =1}^N   ( v_i , v_j ) u_i \otimes u_j.  \]
Where $ S_t^2 ( {\rm Im } \, f_x )  $ denote the subspace of  symmetric tensors in $ {\rm Im} \, f_x \otimes_{F[[t]]} {\rm Im} \, f_x $.
 We remark that though $ {\rm Im} \, f_x $ is a submodule of $M_-$, but the natural map 
  $ S_t^2 (  {\rm Im} \, f_x ) \to    S_t^2 ( M_- )$ is in general {\it not} an embedding.

\begin{theorem}\label{theorem5.1} Two elements $x, y \in M_-\otimes V$
 are in the same $G( F[[t]])$-orbit iff $\, $  $ {\rm Im } \, f_x = {\rm Im } \, f_y $
 and $ T ( x ) = T( y ) $ .
\end{theorem}

\

\noindent We need some preparations for proving the theorem.

\

We recall a special case of Witt's theorem (see \cite{J} ):

\begin{theorem}\label{witt}
If $L_{1},$ $L_{2}\subseteq V[[t]]$ are two primitive submodules of $V[[t]]$ and if 
$\sigma :L_{1}\rightarrow L_{2}$ is an isometry, then $\sigma $ can be
extended to an isometry in $ G ( F[[t]] ) $.
\end{theorem}

\

It is clear that if $x,y$ are in the same $G(F[[t]])$-orbit, then ${\rm Im}\, f_x=
{\rm Im}\, f_y$ and $T(x)=T(y).$  Conversely, if ${\rm Im}\, f_x={\rm Im} \, f_y \deff W$ and 
$T(x)=T(y)$.  Let $e_{1},....,e_{m}$ be a quasi-basis for $W$,
 and $W $ be as ({\ref{5.1.1}).   We may assume that $k_1 \geq k_2 \geq \dots \geq k_m $. Then $S_t^2 ( W)$ has a quasi-basis 
  $  e_{ij} \deff  e_i \otimes e_j + e_j \otimes e_i $ ($ 1 \leq i \leq j \leq m $), and    
 \[  S_t^2 ( W) = \sum_{ 1\leq i \leq j \leq m } F[[t]]/ ( t^{k_{j}} ) e_{ij} .\] 
  By lemma \ref{lemma5.2}, we may write 
 \[  x = e_1 \otimes a_1 + \dots + e_m \otimes a_m  , \, \, \, \, \, \, y = e_1 \otimes b_1 + \dots + e_m \otimes b_m \]
with $\{ a_1 , \dots ,a_m \}$ and $\{ b_1 , \dots , b_m \}$ satisfy  condition (1) in Lemma \ref{lemma5.2}.
Then $ T ( x ) = T( y ) $ implies that 
 \begin{equation}\label{5.8} 
       ( a_i , a_j ) = ( b_i , b_ j ) \, \, \, \, \, {\rm mod} \, \, t^{{\rm min} (k_i , k_j )  } .
\end{equation}

\begin{lemma}\label{lemma5.4} Let $L_1$ and $L_2$ be two primitive submodules of $V[[t]]$
 with bases $a_1 , \dots , a_m $ and $ b_1 , \dots , b_m $. Let $1 \geq k_1 \geq \dots \geq k_m $.
 If (\ref{5.8}) holds, then  the set $ b_1 , \dots , b_m $ can be altered to another set 
 $ {\tilde b}_1 , \dots , {\tilde b}_m $ such that 
  \begin{equation}\label{5.9} 
 {\tilde b}_i =  b_i \, \, \, \, \, \, {\rm mod } \, t^{k_i }  , \, \, \, \, {\rm for } \, 1\leq i \leq m ,
  \end{equation} 
and 
 \begin{equation}\label{5.10} 
  ( {\tilde b}_i , {\tilde b}_i ) = ( a_i , a_j ) \, \, \, \, \, \, {\rm for } \,  1\leq i , j \leq m ,
  \end{equation} 
and the $F[[t]]$-span of $  {\tilde b}_1 , \dots , {\tilde b}_m$  
  is a primitive submodule of $V[[t]]$ with  ${\tilde b}_1 , \dots , {\tilde b}_m$
 as a basis.
\end{lemma}

\

 Suppose the truth of Lemma \ref{lemma5.4}, then Theorem \ref{theorem5.1} can be proved as follows.
 The equation (\ref{5.9}) implies that 
 \[ y = \sum e_i \otimes b_i = \sum e_i \otimes \tilde{b}_i . \] 
 The equation (\ref{5.10}) implies that 
 the map $\sigma :L_{1}\rightarrow L_{2}$ give by $ a_i \mapsto {\tilde b}_i $ is an isometry,
 by Theorem \ref{witt},  $\sigma $ can be extended to $g\in G ( F[[t]])$.
 then 
\[  y \cdot g =  ( \sum e_i \otimes \tilde{b}_i  )\cdot g = \sum e_i \otimes  \sigma (\tilde{b}_i ) = \sum e_i \otimes {a}_i = x .\]  
 It remains to prove Lemma \ref{lemma5.4}.

\

\noindent{\it Proof of Lemma \ref{lemma5.4}.}  We use induction on $m$. For case $m=1$, we first take $c\in V[[t]]$ such that
  $ ( b_1 , c ) = 1 $,  we want to find 
    $ {\tilde b}_1 = b_1 + t^{k_1} h(t) c $
 where $h(t) = h_0 + h_1 t + h_2 t^2 + \dots \in F[[t]] $ such that
  We have 
  \[  (a_1 , a_1 ) =  ( {\tilde b}_1 , {\tilde b}_1 )\]
which is equivalent to  
 \begin{equation}\label{5.11}
   (a_1 , a_1 ) - ( b_1 , b_1 ) = 2 t^{k_1} h( t ) +  t^{2k_1 } h(t)^2 ( c , c )  .
 \end{equation}
Since $ (a_1 , a_1 ) = ( b_1 , b_1 ) \, \, {\rm mod} \, t^{k_1 } $,  we see that 
 \ref{5.11} holds mod $t^{k_1}$ for arbitrary $h(t)$.   Compare the coefficient of $t^{k_1 }$, we solve for 
 $h_0$, after $h_0$, we compare coefficient of $t^{k_1+1}$, we solve $h_1$. It is clear that the similar process can be continued to 
 solve all $h_i$.  Assume the Lemma holds for $m-1$, so we can find $ {\tilde b}_1 , \dots , {\tilde b}_{m-1} $ such that 
  
     \begin{equation}\label{5.9r} 
 {\tilde b}_i \cong b_i \, \, \, \, \, \, {\rm mod } \, t^{k_i }  , \, \, \, \, {\rm for } \, 1\leq i \leq m-1 ,
  \end{equation} 
and 
 \begin{equation}\label{5.10r} 
  ( {\tilde b}_i , {\tilde b}_i ) = ( a_i , a_j ) \, \, \, \, \, \, {\rm for } \,  1\leq i , j \leq m-1 .
  \end{equation} 
 We may assume $b_i = {\tilde b}_i $ for $i=1 , \dots , m-1$.
Since  $(\, )$ on $V[[t]]$ is non-degenerate, we can find 
 $c_1 , c_2 , \dots ,c_m \in V[[t]]$ such that 
\begin{equation}\label{del} 
  ( b_i , c_j ) = \delta_{i, j } .
\end{equation} 
 We want to find $h_1 ( t) , \dots , h_m ( t) \in F[[t]] $ such that 
\[ {\tilde b}_m = b_m + t^{k_m} \left( h_1 ( t ) c_1 + \dots  + h_m ( t)  c_m \right) \]
 satisfies  the following $m$ equations 
      \begin{equation}\label{5.12}   ( {\tilde b}_i , {\tilde b}_m ) = ( a_i , a_m ), \, \, \,  \, \, \, \, i = 1, \dots , m  
     \end{equation} 
Let $h_i ( t) =\sum_{s=0}^{\infty}  h_{i, s} t^s $. 
 The equations (\ref{5.12}) already hold mod $t^{k_m}$. Compare the coefficient of $t^{k_m}$ of (\ref{5.12}), 
 we get a linear system with $m$-variables $h_{1, 0} , \dots , h_{m, 0}$ and $m$ equations, this 
 system has non-zero determinant, thanks to (\ref{del}), we can solve for $h_{1, 0} , \dots , h_{m, 0}$.
 Then we compare coefficient of $t^{k_m+1 }$ of (\ref{5.12}), we get a linear system with $m$ equations
 and $m$ variables $h_{1, 1} , \dots , h_{m, 1}$, and again because of (\ref{del}), the system has a solution.
 It is clear that this process can be continued to solve for all $h_{i, s} $'s.
 \hfill $\Box$

\

\begin{lemma}\label{submersive}  Let $v$ be a place of $F$, $W$ be a $F[[t]]$-submodule of $M_-$, let 
  $W_v = W\otimes F_v, V_v = V\otimes F_v $. 
  The map $T: W_v \otimes V_v \to S_t^2 ( W_v )  $ given by 
 \[ T( \sum_i u_i \otimes v_i ) = \sum_{i, j} ( v_i , v_j ) u_i \otimes u_ j \] 
 is submersive at $x_0 \in W_v \otimes V_v $ iff ${\rm Im } \, f_{x_0} = W_v $
\end{lemma}

\noindent {\it Proof.} For simplicity, we denote $W_v, V_v, F_v$ by $W, V, F$ respectively.
 Let $a_1 , \dots , a_m $ be a quasi-basis of the $F[[t]]$-module ${\rm Im} \, f_{x_0} $. 
    By Lemma  \ref{lemma5.2},  $x_0$ can be written as 
\[  x_0 = \sum_i a_i\otimes b_i = a_1 \otimes b_1 + \dots + a_m \otimes b_m  \] 
where $b_1 , \dots ,b_m$ is a basis of $Span \{ b_1 , \dots ,b_m \}$ and $Span \{ b_1 , \dots ,b_m \}$ is a primitive 
$F[[t]]$-submodule of $V[[t]]$.
  We first find a formula for  the tangent map 
\[ dT_{x_0} :  T_{x_0} =  W\otimes_{F[[t]]} V[[t]] \to  T_{y_0}=  S^2_t ( W ) \]
where $y_0 = T(x_0 )$.
Take a line $  x(\epsilon)=   \sum_i a_i\otimes b_i  + \epsilon \sum_j {u_j \otimes v_j } $ passing through
 $x_0$ in the direction $\sum_j {u_j \otimes v_j } $,
Then 
\[ T ( x(\epsilon) ) = 
\sum ( b_i , b_j) a_i\otimes a_j + \epsilon \sum  (b_i , v_j ) ( a_i \otimes u_j + u_j \otimes a_i ) + 
   \epsilon^2 \sum ( v_i , v_j) u_i\otimes u_j  . \]
So  we have 
\begin{equation}\label{tangent}
 d T_{x_0} (  \sum_j {u_j \otimes v_j } ) =   \sum  (b_i , v_j ) ( a_i \otimes u_j + u_j \otimes a_i ) .
 \end{equation}
 From this formula,  we see that $dT_{x_0}$ is $F[[t]]$-linear.
  If ${\rm Im } \, f_{x_0}\not= W $, then $ {\rm Im} \, dT_{x_0} \subset   {\rm Im } \, f_{x^0 } \otimes W + {\rm Im } \, f_{x^0}\otimes W$,
 $ dT_{x_0}$ is not surjective, i.e. $T$ is not submersive at $x_0$. 
  If  ${\rm Im } \, f_{x^0}= W $, since $( \, )$ on $V[[ t]]$ is non-degenerate, for each $1\leq k \leq m $, we can find 
 $ v \in V[[t]]$ such that $ ( b_i , v) = \delta_{i, k } $, then 
  $ dT_{x_0} (  a_l \otimes v ) = a_l \otimes a_k + a_k \otimes a_l $.
 So   all $ a_l \otimes a_k + a_k \otimes a_l $ are in   ${\rm Im } \, f_{x_0}$. So 
  $dT_{x_0}$ is surjective and $T$ is submersive at $x_0 $.
\hfill $\Box $

 \

 For an $F[[t]]$-submodule $W$ of $M_-$,  $T_W :  W\otimes V \to S_t^2 ( W)$ as in Section 3. For $ i \in  S_t^2 ( W)$,
 We denote by $ U( i)$ the variety of elements in $T_W ^{-1} ( i)$ where $T_W$ is submerssive. 
 Of course the set of $F$-points  $ U(i)_F $ may be empty. By Lemma \ref{submersive}, Theorem \ref{theorem5.1} can be reformulated as

  \begin{theorem}\label{theorem5.2} The $G(F[[t]])$-orbits in $M_-\otimes V$ are in one-to-one correspondence with
    the set of pairs $(W, i )$ with $W\in Gr ( M_- , t )$, $i \in S_t^2 ( W )$ such that $ U(i)_F $ is not empty. 
    The correspondence is the following, for $x \in M_-\otimes V$, its orbit corresponds to the pair $(W, i )$, where
     $ W = {\rm Im} \, f_x $, and $ i = T_W ( x ) $. 
\end{theorem}

\

\

\section{Theta series.}

We continue with the snt-module $M$ with decomposition $ M= M_- \oplus M_+$ and $V$ a finite dimensional vector space
 with a non-degenerate, bilinear symmetric form $(\, ,\, )$ as in Section 3. 
 For each  $\phi \in {\cal S} ( X_{\bf A}) = {\cal S} ( (M_-\otimes V )_{\bf A})$,  the  theta functional $\theta (\phi )$ is 
 defined by (\ref{th}).
Recall $G^q ( {\bf A} [[ t ]] )$ acts on $ {\cal S} ( X_{\bf A})$,  
 the action formula is given as follows.  An element $ g \in G^q ( {\bf A} [[ t ]] )$  has  
  the block decomposition 
\begin{equation*}
\left[ 
\begin{array}{cc}
\alpha _{g} & \beta _{g} \\ 
\gamma _{g} & \delta _{g}%
\end{array}%
\right]
\end{equation*}%
with respect the decomposition 
   \[ (M\otimes V)_{\bf A} = (M_-\otimes V)_{\bf A} \oplus (M_+\otimes V)_{\bf A} .\]
 Since $ G^q ( {\bf A} [[ t ]] )$ preserves $ (M_{\pm} \otimes V)_{\bf A}$, so we have 
 $\beta_g= 0 $ and $\gamma_g = 0 $ .
 Then the action of $g$ on ${ \cal S}  ( ( M_-\otimes V)_{\bf A} )$ is given by 
\begin{equation}\label{4.1}
  ( g \cdot \phi ) ( x ) =  \phi ( x \alpha_g ) =  \phi ( x g )
\end{equation}

If $g\in G^q ( F [[t]] ) $,  it is clear that 
\[   \theta ( g \cdot \phi ) = \theta ( \phi ). \]
So $ \theta ( g \cdot \phi )$ is a continuous function on $G^q ( F [[t]] ) \backslash G^q ( {\bf A} [[t]])$.
 Let  $dg$ be the Haar measure on $ G^q ( {\bf A} [[t]])$  such that the volume of 
 $G^q ( F [[t]] ) \backslash G^q ( {\bf A} [[t]])$ is $1$.

\begin{lemma} \label{lemma6a}
 If $ ( V , ( \, , \, ))$  is anisotropic over $F$ or   ${\rm dim } V  - r > \frac 1 2 {\rm dim} \, M  + 1   $,
 where $r$ is the dimension of a maximal isotropic subspace of $V$,   then 
 the 
integral
\begin{equation}\label{4.5}
    {\rm It} ( \phi  ) \deff \int_{G^q ( F [[t]] ) \backslash G^q ( {\bf A} [[t]]) } \theta (  g \cdot \phi )   dg 
\end{equation}
converges
\end{lemma}

\noindent {\it Proof.}  If $ ( V , ( \, , \, ))$  is anisotropic over $F$, then 
  $G^q ( F [[t]] ) \backslash G^q ( {\bf A} [[t]])$ is compact, so (\ref{4.5}) converges. 
  Let 
  $ U$ be the unipotent radical of $G^q$, then $ G^q = G \ltimes U$. Let $D$ be a compact fundamental domain of 
  $ U ( F ) \backslash U ( {\bf A } )$.  Then 
 \begin{equation}\label{6d}
 (\ref{4.5} ) = \int_{G (F)   \backslash G^q ( A ) } \int_{a\in D } 
    \sum_{r \in M_- \otimes V } \phi ( r a g ) d a    dg
  =   \int_{G (F)   \backslash G^q ( A ) }  \theta ( g \cdot \bar \phi ) dg   , 
 \end{equation}
 where $ \bar \phi ( x ) = \int_{ D} \phi ( x a ) d a  $. Since $D$ is compact, 
 $ \bar \phi \in {\cal S} ( (M_- \otimes V)_{\bf A} ) $. By the convergence criterion in Proposition 8 \cite{W2}, 
 the right hand side of (\ref{6d}) converges under the condition  ${\rm dim } V  - r > \frac 1 2 {\rm dim} \, M  + 1   $.
\hfill $\Box $

We can write the integral (\ref{4.5}) as a sum of orbital integrals. Let ${\cal O}$ be a set of representatives
 of $G(F[[t]])$-orbit in $ M_-\otimes V$. For each $\xi \in {\cal O}$, let $G_{\xi }$ denote  its isotropy subgroup 
 in $G^q(F[[t]])$.  Then 
 the integral (\ref{4.5}) can be written as
\begin{equation}\label{4.7}
\int_{G^q ( F [[t]] ) \backslash G^q ( {\bf A} [[t]]) }   \sum_{\xi \in {\cal O} }
 \sum_{  \tau \in  G_{\xi } \backslash  G^q ( F [[t]])}           \phi (\xi \tau g)dg,
\end{equation}
which can be further written as 

\begin{equation}\label{th3}
 \sum_{\xi \in \mathcal{O}}vol(G_{\xi }\backslash G_{\xi ,{\bf A} } ) 
\int_{G_{\xi ,{\bf A}}\backslash G^{q} ( {\bf A} [[t]])} \phi ( \xi g) dg
\end{equation}%
and we have thereby expressed $ {\rm It } ( \phi )$ as   
  
\begin{equation}\label{orbitint}
 {\rm It } ( \phi ) = \sum_{\xi \in \mathcal{O}}vol(G_{\xi }\backslash G_{\xi ,{\bf A} } ) 
\int_{G_{\xi ,{\bf A}}\backslash G^{q} ( {\bf A} [[t]])} \phi ( \xi g ) dg .
\end{equation}

\

\

\maketitle \section{Siegel-Weil formula}

We assume in this section $M$ and $V$ satisfies the conditions
 that ${\rm dim} \, V > 6 n +2 $, where $n$ is the number of $H_k$'s in the decomposition of $M$ as in (\ref{M3})
 and $V$ satisfies the conditions in Lemma \ref{lemma6a}. 
  By Theorem \ref{conv},
 the $t$-Eisenstein series $\phi \in {\cal S} (  (M_-\otimes V )_{\bf A} )  \mapsto {\rm Et}( \phi ) $ is a tempered distribution on 
 $  (M_-\otimes V )_{\bf A} $.  And for each $ W\in Gr ( M_- , t )$, we have a 
  tempered distribution $\phi \mapsto {\rm Et}_W ( \phi ) $,
  given by (\ref{3ei}). We have 
\[  {\rm Et} = \sum_{W \in Gr (M_- , t ) }  {\rm Et}_W   .\] 
 By Theorem \ref{thm4.7}, 
  \[ {\rm Et}_W = \sum_{ i \in {\rm St}^2 ( W) } \mu_i  .\] 
 We denote $\mu_i $ by ${\rm Et}_{W, i } $.
 Therefore we have 
  \begin{equation}\label{7.2}
 {\rm Et }= \sum_{ W \in Gr ( M_- , t ) } \sum_{ i \in {\rm St}^2 ( W) } \mu_{W, i} .
  \end{equation} 
 Moreover the measure ${\rm Et}_{W, i } = \mu_i $ is the gauge measure
 as described in Theorem \ref{thm4.7}, which implies in particular $\mu_i$ is $0$
 if $ U( i )_{\bf A} $ in empty.

On the other hand,   
\[ \phi \in {\cal S} (  (M_-\otimes V )_{\bf A} ) \mapsto {\rm It} ( \phi )  \]
 given in (\ref{4.5}) is a tempered distribution and it has a decomposition given by (\ref{orbitint}).
 By Theorem \ref{theorem5.2}, each orbit  corresponds uniquely to a pair
 $(W, i )$ where $W \in  Gr( U_- , t)$ and $ i \in {\rm St}^2 ( W)$.
So we may write 
 \[  {\rm It} = \sum_{ W , i } {\rm I}_{ W , i } \]
 where 
\[ {\rm I}_{ W, i } ( \phi ) =  vol(G_{\xi }\backslash G_{\xi ,{\bf A}%
})\int_{G_{\xi ,{\bf A}}\backslash G^{q} ( {\bf A} [[t]])} \phi
(\xi g)dg.  \]
where $(W, i)$ corresponds to the orbit containing $\xi $, i.e., 
 ${\rm Im } \, f_{\xi } = W $ and $T_W ( \xi ) = i $. 
 If $U(i)$ is empty, $(W, i)$ doesn't corresponds to any orbit, in this case we define 
  \[ {\rm I}_{ W, i }\deff 0 .\]
We shall prove that ${\rm Et} ={\rm  It }$, and actually we shall prove more:
 $ {\rm Et}_{W , i } = {\rm It}_{W, i}$ for any pair $(W, i )$
. We use the induction on ${\rm dim } \, M$.
 The case ${\dim } \,  M =2 $ is the classical result  in \cite{W2}.

 Our proof is entirely parallel to that of \cite{W2}. To start with, we introduce some notations.
 Let $\pi : \widehat{Sp}_{2N} ( {\bf A}) \to  Sp_{2N} ( {\bf A})$ denote the double cover (recall
 $2N ={\rm dim } \, M {\rm dim} \, V $). Let $  \widehat{Sp} ( M , t )_{\bf A} $ denote
  $\pi^{-1}  ( {Sp} ( M , t )_{\bf A} )$.  We let 
 \[  M ( k ) = \oplus_{k_i = k } H_{k_i } ,\]
 so  
\[ M =  M( l_1 ) \oplus \dots \oplus M( l_s ) , \, \, \, \, \, l_1 > \dots >l_s . \]
 Recall Corollary 2.6, $Sp( M , t ) = N \ltimes H$,  where $N$ is the unipotent radical,
  and 
 \[ H  = \Pi_{i=1}^s Sp_{2r_i } ( F ) ,\]
 where $r_i$ is the number $H_{k_i}$'s in the decomposition of $M_{l_i}$. 
Since $M$ has decomposition (\ref{M3}), $M_-$ has decomposition
\begin{equation}\label{M-}
  M_- = F [ t] / ( t^{k_1}) \oplus \dots \oplus   F [ t] / ( t^{k_n}) .\end{equation}
 Let $T$ be a maximal torus of $H$,  we may take $T$ such that $T$  preserves $M_{\pm}$ 
 and preserves each component in (\ref{M-}).
 Then  $T = {\rm G}_{\rm m}^n$, where the $i$-th $G_m$ acts on the $i$-th component in  (\ref{M-}).
  We have $T_{\bf A} =  I_F^n $, where $I_F$ denotes the idele group of $F$.
  Let $T_{\bf A}'$ be the subset of $T_{\bf A}$ formed by 
\[  T_{\bf A}' = \{ ( t_1 , \dots , t_n ) \, | \, |t_1|_{\bf A} \geq \dots \geq  |t_{r_1}|_{\bf A} \geq 1,
  \dots ,  |t_{n-r_s+1 }1|_{\bf A} \geq \dots \geq  |t_{n}|_{\bf A} \geq 1
   \}. \]
 Since $N$ is unipotent, the space $ N( F )  \backslash N ( {\bf A} ) $ is compact. By the reduction theory for the
 semi-simple group $H$, there exists a compact  $C \subset  Sp( M , t )_{\bf A}$ such that  
   \begin{equation}\label{red}
     Sp( M , t )_{\bf A} =  Sp ( M , t)_F T_{\bf A}' C .
  \end{equation} 
 Since $T$ preserves $M_+$, we may regard $T_{\bf A}$ as a subgroup of $\widehat {  Sp} ( M , t )_{\bf A}$,
 the decomposition (\ref{red}) implies that 
\begin{equation}\label{red2}
    {\widehat Sp}( M , t )_{\bf A} =  Sp ( M , t)_F T_{\bf A}' C .
  \end{equation} 
for some compact subset $C\subset {\widehat Sp}( M , t )_{\bf A}$.

As in \cite{W2}, for each $\tau \in {\Bbb R}_{ > 0 }$, we let
 $a_{\tau }\in I_F$  denote the idele such that $ (a_{\tau })_v  = \tau $ for each infinite place $v$
  and $ (a_{\tau })_v  = 1 $ for each finite place.  We let
 $ \Theta ( T)$ denote the set of all $ ( a_{\tau_1} , \dots , a_{\tau_n } )$, and set 
  $  \Theta ( T)' = \Theta ( T) \cap T_{\bf A}' $.

\begin{lemma}\label{lemma7.1}
 If $\hat E$ is a positive tempered measure on $ X_{\bf A} =  (M_-\otimes V)_{\bf A}$,  and is a sum of positive measures $\hat \mu_i$ 
supported on $ U(i)_{\bf A}$ {\rm (}$i\in S_t^2 ( M_-)$ {\rm)}, and  is
 $T_F$-invariant,  and  
 there is a place $v$ of $F$ and a subgroup $G_v' $ of $G^q ( F_v [[ t]])$ that acts transitively on $U(i)_v$ 
   such that $\hat E$ is invariant under $G_v'$. Then the function 
 $ S \mapsto \hat E ( S \phi )$ is bounded on $T'_{\bf A}$, uniformly for $\phi$ in a compact subset 
 in ${\cal S} ( X_{\bf A} )$.
\end{lemma}

This lemma is a generalization of Lemma 23 in \cite{W2}.  Our proof below closely follows that of \cite{W2}. 

\noindent {\it Proof.} Let $e_1 , \dots , e_n$ be a quasi-basis of $M_-$:
  \[ M_- = F [t] / ( t^{k_1 } ) e_1 \oplus \dots \oplus  F [t] / ( t^{k_n } ) e_n. \]
 Then $e_{i} \otimes e_j +e_{j} \otimes e_i $ is a quasi-basis of $ {\rm St}^2 ( M_-)$.
 For each $\alpha \in \{ 0 , 1, \dots , n \}$, let $ S_t^2 ( M_-)^{ ( \alpha ) } $ be the 
 set that consists of elements 
\[  \sum_{ i , j > \alpha  } k_{ij }(  e_i\otimes e_j +  e_{j } \otimes e_ i ) \] 
such that 
\[  k_{\alpha+1 , j } ( e_{\alpha + 1 } \otimes e_j +  e_{j } \otimes e_{\alpha  + 1 } ) \ne 0 \] 
for at least one $j \geq {\alpha + 1 } $. 
 We set
\[ S_t^2 ( M_-)^{ ( n ) } = \{ 0 \} \]
by convention.
  It is clear that $ S_t^2 ( M_-)$ is a disjoint union of $ S_t^2 ( M_-)^{ ( \alpha ) }$.
 Let $\hat E_{\alpha }$ be the sum of $\hat \mu_i $ for $i \in   S_t^2 ( M_-)^{ ( \alpha ) }$.
 It is clear that 
 \[ \hat E = \hat E_{0 } +  \hat E_{1 } + \dots + \hat E_{n } \]
and $\hat E_{\alpha }$ satisfies all the conditions in the lemma. it is enough to prove 
 the result for each $\hat E_{\alpha } $.

Now we fix $0 \leq \alpha \leq n $.
Let $q$ be the constant which is $1$ if $v$ is infinite and is equal to the cardinality of the residue
 field if $v$ is a finite place. As in \cite{W2}, there is a compact subset $ C \subset I_F$ (where $I_F$ is the idele group of $F$), such that 
 every $ t\in I_F$ with $    1\leq  | t | \leq S $ (where $S$ is a fixed constant )
  can be written as $r c $ with $r \in F$, $c \in C$.
We denote $C^n$ the subset of 
  $T_{\bf A}$ formed by elements $( c_1 , \dots , c_n)$ with all $c_i\in C$,
 and let 
 \[ \Theta'_{\alpha } = \{ ( a_{\tau_1 } , \dots ,a_{\tau_n} ) \, | \, \tau_1 = \dots = \tau_{\alpha +1 } \geq \dots \geq \tau_n \geq 1 \}.\]

We will apply Lemma 6 of \cite{W2} to the space
 $X_F $. We consider $X_F$ as a product space $\Pi_{i=1}^{n-\alpha } X_F^{(i)} $, where
 $X_F^{(1)}$ is the $F[[t]]$-submodule generated by 
 $ e_1 \otimes V , \dots , e_{\alpha +1 } \otimes V$, and for $n-\alpha \geq i\geq 2$, 
  $   X_F^{(i)}$ is  $ F[[ t]] e_{ i + \alpha } $,
  $Y_F \subset M_{-} \otimes_{ F [[ t]] } M_-$  
 is the submodule over $F[[t]]$ spanned by $ e_{i} \otimes e_j  $ with $j \leq \alpha +1 $ .    
   And $p : X_F \to Y_F $ is given by 
 \[  p ( \sum_{i=1}^{j=1}  e_i \otimes v_ i  ) = \sum_{i=1}^n \sum_{j=1}^{\alpha+1} ( v_i , v_j) e_i \otimes e_j   .\]
 Apply Lemma 6 of \cite{W2}, there is $\phi_0 \in S( X_{\bf A} )$ such that 
  \begin{equation}\label{7.5}
  | ( ( \theta c ) \cdot \phi ) (x )  | \leq \phi_0 ( x )  
 \end{equation}
for all $x\in X_{\bf A}$ with $T ( x ) \in S_t ( M_-)^{(\alpha)}$, all $\theta \in \Theta'_{\alpha }$, $c\in C^n $,
 and $\phi \in C_0$.

 It can be proved (\cite{W2} page 71) that 
 each element 
 $ t = ( t_1 , \dots , t_n )\in T'_{\bf A} $ can be written as
 \[ t = r y \theta c  \]
where $ r \in T_F $, $y = ( y_1 , \dots , y_n ) \in T_v $
 with all $|y_i|_v \geq 1 $ and $y_{\alpha +1 } = \dots = y_n = 1 $, $\theta \in \Theta'_{\alpha }$ and $c\in C^n$. 
 Since $\hat E_{\alpha }$ is invariant under $T_F$, we have 
\[ |  \hat E_{\alpha } ( t \cdot \phi ) | \leq \hat E_{\alpha} ( y \cdot \phi_0 ) . \] 
We may assume that $\phi_0 = \phi_v \phi' $ where $\phi_v \in {\cal S} ( X_v )$, 
 $\phi' \in   {\cal S} (\Pi'_{w\ne v }  X_w )$. By Lemma 22 \cite{W2}, we have 
\[ \hat E_{\alpha} ( y \cdot \phi_0 ) = \sum_{i \in St^2 ( M_-)^{(\alpha )}} c_i ( \phi ' ) 
 \int_{U(i)_v } y \cdot \phi_v | \theta_i|_v .\]
From this, we obtain that 
  
\[ \hat E_{\alpha} ( y \cdot \phi_0 ) = \Pi_{i=1}^n |y_i|_v^{ 	 k_1 + \dots + k_{i-1} + ( n - m/2 - i +2  )k_i }
  \hat E_{\alpha} (  \phi_0 ) \]
since $k_i \leq k_j $ for $i \leq j$, and $m > 6 n +2 $, we see that 
 the exponent of $|y_i|_v$ is $\leq 0$, and since $|y_i|_v \geq 1 $, so have 

\[  \hat E_{\alpha} ( y \cdot \phi_0 ) \leq \hat E_{\alpha } ( \phi_0 ) .\]
This proves the lemma. \hfill $\Box $

\

The following theorem is a generalization of Theorem 4 in \cite{W2}.

\begin{theorem}\label{unique} Suppose $ {\rm dim } \, V  > 6 n + 2$. If there is a place $v$ of $F$ such that
 $ U ( 0 )_v $ is not empty, and a subgroup $G_v'$ of $G^q ( F_v [[ t]] )$ acts transitively on $U(i)_v$ for 
  every $ i\in S_t^2 ( M_-)$. And if $E'$ is a positive tempered measure on $X_{\bf A} $ invariant under $Sp ( M , t) $ and 
  $G_v'$ and $ E' -{\rm Et } $ is supported on the union of  $U(i)_{\bf A}$ for $ i \in S_t^2 ( M_- )$ .
 Then $ E' = {\rm Et} $. 
\end{theorem} 

 \noindent {\it Sketch of Proof. }   Using Lemma \ref{lemma7.1} and \ref{red2}, 
  it is easy to see  that
  the function $ \widehat{Sp} ( M , t )_{\bf A} \to {\Bbb C}$ given by 
 $ S \mapsto  ( E' - {\rm Et } ) ( S \phi )$ is bounded on $\widehat{Sp} ( M , t )_{\bf A}$ uniformly for 
 $\phi $ in every compact subset of ${\cal S} ( X_{\bf A} )$.  The remainder of the proof is similar to that of Theorem 4 in \cite{W2}.
\hfill $\Box $

\

Now we can prove the main theorem of this work:

\begin{theorem}\label{main} Suppose $ {\rm dim } \, V  > 6 n + 2$ and $V$ is anisotropic or   ${\rm dim } V  - r > \frac 1 2 {\rm dim} \, M  + 1   $,
 where $r$ is the dimension of a maximal isotropic subspace of $V$, then 
  \[ {\rm Et } = {\rm It }.  \]
\end{theorem} 

 The condition on $V$ in the Theorem is for the convergence of ${\rm It} $ (see \ref{lemma6a}).  
The proof uses the induction on ${\rm dim } \, M $.  The induction assumption 
 implies that $ E' =  {\rm It } $ satisfies the conditions of Theorem  \ref{unique}, therefore
  $  {\rm Et } = {\rm It }$.

\

\

\maketitle \section{Corollaries of Siegel-Weil formula  }

\

In this section we prove a slightly more general form of the Siegel-Weil formula (Theorem \ref{thm8.1}) for snt-modules 
 that will be used in part II.  
 Let $M,  V$ be as in Section 3. Let 
 \begin{equation}\label{8.1}
 M = M_- \oplus M_+ 
\end{equation}
 be a direct sum such that  $M_+ \in Gr ( M , t )$, $M_-\in Gr ( M)$ but not necessarily in $Gr ( M , t )$.  
 Let $ X = M_- \otimes V $.
 The space 
  $L^2 ( X_{\bf A}) $ is a representation of the
  metaplectic group  ${\widehat Sp}_{2N} ( {\bf A})$ ($2N = {\rm dim} M  {\rm dim } V $) with the usual  theta functional 
\[ \theta :   {\cal S} (  X_{\bf A}) \to {\Bbb C} , \, \, \, \phi \mapsto \theta ( \phi ) = \sum_{r\in X } \phi( r) .\] 
  Recall  the Eisenstein series  (\ref{eisen1}), (\ref{eq1s}) for $\phi \in {\cal S} ( X_{\bf A} )$ is given by 
\[ {\rm E } ( \phi ) = \sum_{ U \in Gr ( M )} E(\phi , U) =
  \sum_{U \in Gr ( M) }  \int_{ (\pi_- ( U )  \otimes V )_{\bf A} }   \psi (  \frac 12 \la x , \rho x \ra )     \phi ( x ) dx  .   \]
The $t$-Eisenstein series defined in (\ref{Et}) is a subseries given by  
 \begin{equation}\label{8t}
{\rm  Et} ( \phi )  = \sum_{ U \in Gr ( X , t )} E(\phi , U)  .    
 \end{equation}

 \noindent Though $ M_-$ is not an $F[[t]]$-submodule of $M$,  it has $F[[t]]$-module structure via the isomorphism 
 $M_- = M / M_+$.  In the case that $M_-$ is an $F[[t]]$-submodule, the two $F[[t]]$-module structures on $U_-$ clearly coincide.   
 The projection map 
  \[ \pi_- : M \to  M / M_+ = M_- \] 
 is an $F[[t]]$-module homomorphism. 
 For each $ U \in Gr ( M , t ) $, $\pi_- ( U ) $ is an $F[[t]]$-submodule of $M_-$. 
 Denote $Gr ( M_- , t )$ the set of $F[[t]]$-submodules of $M_-$, so we have a map 
 \begin{equation}\label{3.2}  P :  Gr( M , t ) \to Gr ( M_- , t ) :  \, \, \, U \mapsto \pi_- ( U ) .\end{equation}

For $W\in Gr ( M_- , t ) $, we set 
\[  {\rm Et}_W  ( \phi   ) = \sum_{ U \in Gr ( M, t) :  \pi_- ( U ) = W  } E ( \phi , U ) . \]

We recall the action formula of $G^q ( {\bf A} [[ t ]] )$  on $ {\cal S} ( (M_-\otimes V )_{\bf A})$.  
  An element $ g \in G^q ( {\bf A} [[ t ]] )$  has  
  the block decomposition 
\begin{equation*}
\left[ 
\begin{array}{cc}
\alpha _{g} & \beta _{g} \\ 
\gamma _{g} & \delta _{g}%
\end{array}%
\right]
\end{equation*}%
with respect the decomposition 
   \[  X_{\bf A} = (M_-\otimes V)_{\bf A} \oplus (M_+\otimes V)_{\bf A} .\]
 Since $ G^q ( {\bf A} [[ t ]] )$ preserves $ (M_+\otimes V)_{\bf A}$, so we have $\gamma_g = 0 $.
  Since $ M_-$ is not an $F[[t]]$-submodule in general,  $\beta_g$ may not be $0$.
 Then the action of $g$ on ${ \cal S}  ( ( X\otimes V)_{\bf A} )$ is given by 
\begin{equation}\label{8a}
  ( g \cdot \phi ) ( x ) = \psi ( \frac 12 \la x \alpha_g , x \beta_g \ra ) \phi ( x \alpha_g ) 
\end{equation}

If $g\in G^q ( F [[t]] ) $ and $\xi \in V$, then
\begin{equation*}
\frac{1}{2} \la \xi \alpha _{g}, \xi \beta _{g} \ra \in F,
\end{equation*}%
and so $\psi (\frac{1}{2} \la \xi \alpha _{g}, \xi \beta _{g} \ra )=1$ and%
\begin{equation*}
  ( g  \cdot \phi ) ( \xi )=\phi ( \xi \alpha _{g}).
\end{equation*}
    Therefore,
 we have
\begin{equation}
\theta  (  g\cdot \phi  )= \theta  ( \phi  ) , \, \, \, \, \, {\rm for} \, \,  
 g \in G^q ( F[[t]] ).  
\end{equation}
 The function $ \theta ( g\cdot \phi ) $ (as a function on $G^q ( {\bf A} [[t]] )$
  is actually a function on $ G^q ( F [[ t ]] )  \backslash G^q ( {\bf A} [[ t ]] )$.
 Assume $( V , ( \, , \, ))$ satisfies the conditions in Lemma \ref{lemma6a}, then 
 we can form the convergent
integral%
\begin{equation}\label{8b}
    {\rm It}  ( \phi  ) \deff \int_{G^q ( F [[t]] ) \backslash G^q ( {\bf A} [[t]]) } \theta (  g \cdot \phi ) dg ,   
\end{equation}
the convergence can be proved in the same way as Lemma \ref{lemma6a}.

We wish to write the integral (\ref{8b}) as a sum of orbital integrals. 
We introduce the set
\begin{equation*}
\Omega =S^{1}\times (M_{-}\otimes  V)_{\bf A },
\end{equation*}%
on which $G^{q}( {\bf A} [[t]] )$ acts as 
\begin{equation}\label{8c}
(s,x)\cdot g=(s\psi (\frac{1}{2} \la x\alpha _{g},x\beta _{g} \ra ),x\alpha
_{g}).   \end{equation}
One can then check directly that (\ref{8c}) does define a group action; i.e., that
\begin{equation*}
((s,x)g_{1})g_{2}=(s,x)g_{1}g_{2},\; \, \, \, \, g_{1},g_{2}\in G^{q} ( {\bf A} [[t]]).
\end{equation*}%
On the other hand, if
\begin{equation*}
\phi \in {\cal S}((U_{-}\otimes V)_{\bf A}),
\end{equation*}%
we can extend $\phi $ to $\Omega $ by
\begin{equation*}
\varphi (s,x)=s\varphi (x),\; \, \, \, \, x\in (U_{-}\otimes V)_{{\bf A}}, \, \, \, \, s\in
S^{1}.
\end{equation*}
And note that we can rewrite the integral (\ref{8b}) as
\begin{equation}\label{8b2}
\int_{G^q ( F [[t]] ) \backslash G^q ( {\bf A} [[t]]) }  \left( \sum_{\xi \in U_{-}\otimes V}  \phi ((1,\xi )\cdot g) \right) dg.
\end{equation}
We note that the subset
\begin{equation*}
(1,U_{-}\otimes  V)\subset \Omega
\end{equation*}
is invariant under $G^{q}( F [[t]]) .$  We let $\mathcal{O}\subset
(1,U_{-}\otimes V)$ be a family of orbit representatives for the action
of $G^{q}( F [[t]] ) $.  Of course $\mathcal{O}$ can be identified with a subset
of $U_{-}\otimes _{F}V$ $((1,\xi )\mapsto \xi ,$ $\xi \in U_{-}\otimes V),$ 
and as such, it is a family of orbit representatives for the action%
\begin{equation*}
\xi \longmapsto \xi \alpha _{g},\; \, \, \, \, g\in G^{q} ( F [[t]] ) 
\end{equation*}%
of $G^{q}( F [[t]] ) $ on $U_{-}\otimes V$. 
We
 may
 rewrite
 (\ref{8b2}) as
\begin{equation}\label{8b3}
\int_{ G^q ( F [[t]] ) \backslash G^q ( {\bf A} [[t]])   }\sum_{\xi \in \mathcal{O}}\sum_{\tau \in G_{\xi }\backslash
G^{q}(F[[t]])}\varphi ((1,\xi )\tau  g)dg. 
\end{equation}%
Using the fact that $ G_{\xi ,{\bf A}}$ is unimodular, we can prove  
 (\ref{8b3}) is equal to 
\begin{equation}\label{8b4}
 {\rm It} ( \phi ) 
=\sum_{\xi \in \mathcal{O}}vol(G_{\xi }\backslash G_{\xi ,{\bf A}%
})\int_{G_{\xi ,{\bf A}}\backslash G^{q} ( {\bf A} [[t]])} \phi
((1,\xi )g_{2})dg_{2},
\end{equation}%
by a similar argument as in \cite{W2} (page 16). 
We have thereby expressed the ${\rm It} (\phi ) $  as a sum of
orbital integrals.  For each orbit $ {\cal O}$, we have a
 corresponding $W \in Gr ( M_- , t )$ (see Theorem \ref{theorem5.2}).  Let $ {\rm It }_W ( \phi )$ be the subseries
 of (\ref{8b4}) that is over all $ \xi $ corresponding to $W$. Then we have

\[  {\rm It } ( \phi ) = \sum_{ W \in Gr ( M_- , t ) } {\rm It}_W  ( \phi ) . \]

\begin{theorem}\label{thm8.1} Suppose $ {\rm dim } \, V  > 6 n + 2$ and $(V, (\, , \,))$ satisfies the 
 condition in Lemma \ref{lemma6a}.  
  Let $ M = M_- \oplus M_+ $ be a decomposition such that 
  $M_-\in Gr ( M)$ and $ M_+ \in Gr ( M , t)$. Put $ X = M_- \otimes V$. 
 We can define tempered distributions ${\rm Et}$ and ${\rm It }$ on $X_{\bf A} $ as in
 (\ref{8t}) and   (\ref{8b}). Then 
 \begin{equation}\label{8e1} {\rm Et } = {\rm It }.  \end{equation}
And for each $W\in Gr ( M_- , t )$, we have 
  \begin{equation}\label{8e2} {\rm Et }_W = {\rm It }_W .\end{equation} 
 \end{theorem} 

Notice that we didn't assume $M_-\in Gr ( M , t )$, but this theorem can be reduced to Theorem \ref{main}. 
 Let 
\[  M= \bar M_- \oplus \bar M_+ \]
 be a decomposition of $t$-Lagrangian subspaces.  Both $   L^2 (  ( M_- \otimes V )_{\bf A} )  $ and
  $ L^2 (  (\bar M_- \otimes V )_{\bf A} )$ are models of the Weil representation of $ \widehat{Sp}_{2N} ( {\bf A})$.
 We shall define  an intertwining operator 
  \[   T :    L^2 (  ( M_- \otimes V )_{\bf A} ) \to  L^2 (  (\bar M_- \otimes V )_{\bf A}  )  \]
 such that $T$ sends ${\cal S} (  ( M_- \otimes V )_{\bf A}  ) $ to $ { \cal S}  (   (\bar M_- \otimes V )_{\bf A}  )$, and
\[  {\rm Et } (  \phi )  =  {\rm Et } ( T \phi ) , \, \, \, \, \, \, \, {\rm It } (  \phi )  =  {\rm It } ( T \phi ) .\]
   By Theorem \ref{main},  $ {\rm Et } ( T \phi ) =  {\rm It } ( T \phi )$, so we have ,
  $ {\rm Et } (  \phi ) =  {\rm It } (  \phi )$.  The more detailed proof follows.

\noindent {\it Proof.}  Recall for the symplectic space $  M \otimes V $, 
  we have the associated Heisenberg group 
 \[  H  =   (M\otimes V)_{\bf A}  \times S^1 \]
with the product given by 
 \[   ( a_1 , s_1 ) ( a_2 , s_2 ) = ( a_1 + a_2 ,  s_1 s_2 \psi ( \frac 12  \la a_1 , a_2 \ra ) .\]
Since $M_+$ is a Lagrangian subspace,   $ A \deff  ( M_+ \otimes V )_{\bf A} \times S^1 $ is a maximal 
  abelian subgroup of $H$. And 
 \[ \chi :   ( M_+ \otimes V )_{\bf A} \times S^1  \to S^1 , \, \, \, \, \, 
     ( a , s) \mapsto s \]
 is a $1$-dimensional representation. The induced representation
\[   Ind_{A }^{H } \chi \]
 consists of functions $f$  on $ H$ such that $ f ( a x ) =\chi ( a ) f (x ) $ 
 for all $a \in A$ and the restriction $f$ on $(M_-\otimes V)_{\bf A}$ is in $L^2 ((M_-\otimes V)_{\bf A})   $.
  Let $ M = \bar M_- \oplus \bar M_+$ be a direct sum of $t$-Lagrangian subspaces.
  Similarly, for the maximal abelian subgroup $ ( \bar M_- \otimes V )_{\bf A} \times S^1$, 
   and the character 
 \[  \bar \chi :  ( \bar M_- \otimes V )_{\bf A} \times S^1 \to S^1 , \, \, \, \, \,  ( \bar a , s) \mapsto s \]
 we  have the induced representation 
 \[  Ind_{\bar A }^{H } \bar \chi .\]
The representation space $  Ind_{A }^{H } \chi$ ( $Ind_{\bar A }^{H } \bar \chi$, resp.) 
  can be identified with $ L^2 ( ( M_- \otimes V)_{\bf A} )$ ( $ L^2 ( (\bar M_- \otimes V)_{\bf A} )$, resp.)
  by the restricting a function on $H$ to $( \bar M_- \otimes V )_{\bf A}$ ($( \bar M_- \otimes V )_{\bf A} $, reps.).
  The space of smooth vectors are ${\cal S} (  M_- \otimes V)_{\bf A})$ and  ${\cal S} ( (\bar M_- \otimes V)_{\bf A}) $, 
 respectively. 
Let $  h  \in Sp  ( M )$ be an element such that 
 \[    M_- h  =  \bar  M_- ,   \, \, \, \, \, \, \, \, M_+ h = \bar M_+  .\]
We define 
  \[ T:   {\cal S} (  M_- \otimes V)_{\bf A})  \to {\cal S} ( (\bar M_- \otimes V)_{\bf A}) \]
 by  \[    ( T f ) ( x ) =   ( h  \cdot f ) (  x h^{-1}  )  .  \]  
 It is easy to check that $T$ is an isomorphism of $H$-representations.
There are two actions of $\widehat{{Sp}_{2N}} ( {\bf A} )$ on  
  ${\cal S} ( (\bar M_- \otimes V)_{\bf A})$
(or on $L^2 (   {\bar M }_-\otimes V)_{\bf A} )$): the Weil representation action which we denote by 
 $\pi  ( g ) $,   and $ \pi' ( g ) = T  \pi ( g )  T^{-1}$ (where $\pi ( g )$ denote the Weil representation action 
 on  ${\cal S} ( { M }_-\otimes V)_{\bf A} )$.  
   They both satisfy that,  for every $\alpha \in H:$

 \[  \pi ( g^{-1}  ) \pi ( \alpha )  \pi ( g  ) =   \pi ( \alpha \cdot g )  \]
\[  \pi' ( g^{-1}  ) \pi ( \alpha )  \pi' ( g  ) =   \pi ( \alpha \cdot g )  \]
So $\pi ' ( g ) = c_g \pi ( g )  $ for some scalar $c_g$.   It is clear that $ c_{g_1 g_2 } = c_{g_1 } c_{g_2} $.
  Since $ Sp_{2N} ( F)$ is a perfect group, we have $c_g= 1 $ for all $ g\in  Sp_{2N} ( F)$.
    This implies that $c_g =1 $ for all $g \in \widehat{{Sp}_{2N}} ( {\bf A} ) $. 
It then follows that $T$ is an isomorphism of the Weil representations.
    We define as usual the theta function functionals 
  \[ \theta :  {\cal S} ((  M_- \otimes V)_{\bf A})  \to {\Bbb C} , \, \, \, \, \, \, \,  \bar \theta  : {\cal S} ( (\bar M_- \otimes V)_{\bf A} ) 
  \to {\Bbb C} . \]
  We have  
\[  \bar \theta ( T f ) = \sum_{ r \in \bar M_- }  (T f ) ( r ) 
  =  \sum_{ r \in \bar M_- } (  h \cdot f ) ( r h^{-1} ) =    \sum_{ r \in  M_- } (  h \cdot f ) ( r  )  
  = \theta ( h \cdot f ) = \theta ( f ) .\] 
And 
 \begin{eqnarray}
   {\rm It} (  T f )   & = & \int_{  G^q ( F [[ t]]\backslash  G^q ( {\bf A} [[ t]] }   \theta ( g \cdot T f ) d g     \nonumber  \\
    & = & \int_{  G^q ( F [[ t]]\backslash  G^q ( {\bf A} [[ t]] }  \bar \theta ( T (g \cdot f) ) d g      \nonumber  \\
 & = & \int_{  G^q ( F [[ t]]\backslash  G^q ( {\bf A} [[ t]] }  \theta (  (g \cdot f) ) d g       \nonumber   \\
 & = &  {\rm It} ( f ).      \nonumber 
 \end{eqnarray}

      Let $P \subset Sp( M) $ be the parabolic subgroup that consists of all $ g \in   Sp( M)$ such that $ M_+  g = M_+$,
 similarly,  let $\bar P \subset Sp( M) $ be the parabolic subgroup that consists of all $ g \in   Sp( M)$ such that $ \bar M_+  g =\bar M_+$.
  And let $ S$ denote the set of all $ g \in  P \backslash Sp( M)$ such that $ M_+ g$ is a $t$-Lagrangian subspace.
 Similarly let $\bar  S$ denote the set of all $ g \in \bar P \backslash Sp( M)$ such that $ \bar M_+ g$ is a $t$-Lagrangian subspace.
   Since $   M_+ h = \bar M_+ $, the map $ L_h :  \bar S \to  S,    g \mapsto hg $ is a bijection.

\begin{eqnarray}   {\rm Et } ( T f ) &=& \sum_{   g \in \bar S }  ( g \cdot T f ) ( 0 )      \nonumber  \\
    &=& \sum _{   g \in \bar S }  ( T g   f ) ( 0 )       \nonumber    \\  
     &=&      \sum _{   g \in \bar S }  ( h g   f ) ( 0 )    \nonumber    \\
     &=&  \sum _{   g \in  S }  (g   f ) ( 0 )       \nonumber   \\
     &=&  {\rm Et } (  f ). \nonumber
    \end{eqnarray}
 This proves $ {\rm Et} ( \phi ) = {\rm It } ( \phi ) $.  To prove $ {\rm Et}_W ( \phi ) = {\rm It }_W ( \phi )$
 for every $ W \in Gr ( M_- , t ) $,  we first note for every
 $W\in Gr ( M_- , t )$, $ W + M_+ $ is an $F[[t]]$-submodule of $ M$. 
 Let 
  \[ W^{\bot} = \{  v \in M_+ \, | \, \la v , W \ra = 0 \} .\]
 It is easy to see that $ W^{\bot }$ is the radical of the restriction of $\la \, , \, \ra $ on $W + M_+ $.
  Therefore 
\[  M_W \deff  (W + M_+ ) /  W^{\bot }  \]
 is an snt-module. And 
  \[  M_W =  W \oplus M_+ / W^{\bot }  \]
 is a decomposition into Lagrangian subspaces, and $   M_+ / W^{\bot } \in Gr ( M_W , t )$.
 We consider $ {\cal S} ( (W \otimes V )_{\bf A}) $,  
and apply (\ref{8e1}) to  the pair $( Sp ( M_W , t ) , G )$, we have 
\[   \sum_{ W' \in Gr (M_- , t ) : W' \subset W }  {\rm Et}_{W'} ( \phi ) = 
    \sum_{ W' \in Gr (M_- , t ) : W' \subset W }  {\rm It}_{W'} ( \phi ) .  \]
 The above equality holds for all $ W \in Gr ( M_- , t)$, which 
 implies $ {\rm Et}_W ( \phi ) 
     =  {\rm It}_W ( \phi ) $. \hfill $\Box $

Dept of Math, Yale University, New Haven, CT 06520-8283.   hgarland@math.yale.edu 

Dept of Math, Hong Kong University of Science and Tecgnology.    mazhu@ust.hk


\begin{thebibliography}{aa}




\bibitem{Garland1}
H. Garland,  Certain Eisenstein series on loop groups: convergence and the constant term,  
 {\em   Algebraic Groups and Arithmetic}  Tata Inst. Fund. Res. Mumbai (2004), 275-319.

\bibitem{Garland2}
H. Garland, Absolute convergence of Eisenstein series on loop groups,
 {\em Duke Math J.} Vol. 135, No. 2 (2006), 203-260.

\bibitem{GZ}
H. Garland, Y.Zhu, On the Siegel-Weil theorem for loop groups (II), preprint, 2008. 

 \bibitem{J} D.G. James, On Witt's theorem for unimodular,
quadratic forms II, {\em  Pacific Jour. of Math}. 33, 1970.

\bibitem{KR} S.S. Kudla, S. Rallis,  On the Weil-Siegel formula,
 {\em  J. reine angew. Math.  387} (1988), 1-68.


\bibitem{W1}
A. Weil, Sur certaines groupes d'oprateurs unitaries, {\em Acta Math. 111} (1964), 143-211.

\bibitem{W2}
A. Weil, Sur la formule de Siegel dans la theorie des groupes calssiques,  {\em Acta Math. 113} (1965), 1-88.


\bibitem{Z}
Y. Zhu,  Theta functions ans Weil representations of loop symplectic groups, {\em  Duke Math. J}, vol 143 (2208), no. 1, 17-39.


\end{thebibliography}
\end{document}